\documentclass[12pt,onecolumn,journal,letterpaper]{IEEEtran}
                               
\usepackage{soul}
\usepackage{mathtools}
\usepackage{bm}
\usepackage{algorithmicx}
\usepackage{algorithm,setspace}
\usepackage{algpseudocode}
\usepackage{graphicx}
\usepackage{epstopdf}
\usepackage{amsfonts}
\usepackage{amsmath}

\usepackage{mathtools}
\usepackage{epsfig} 
\usepackage[caption=false]{subfig} %

\usepackage{xcolor}
\usepackage{amssymb}

\usepackage{dsfont} %

\newtheorem{theorem}{Theorem}[section]

\newtheorem{lemma}[theorem]{Lemma}

\newtheorem{assumption}[theorem]{Assumption}

\newtheorem{definition}[theorem]{Definition}
\newtheorem{remark}[theorem]{Remark}

\newcommand\oprocendsymbol{\hbox{$\square$}}
\newcommand\oprocend{\relax\ifmmode\else\unskip\hfill\fi\oprocendsymbol}

\usepackage{balance}
\usepackage[printonlyused]{acronym}

\acrodef{MILP}{\emph{Mixed-Integer Linear Program}} %
\acrodef{MILPs}{\emph{Mixed-Integer Linear Programs}} %
\acrodef{MIG}{\emph{Mixed-Integer Gomory}}%
\def \DiMILP/{\texttt{DiCUT-MILP}} %

\usepackage{soul}

\newcommand{\varz}{z}
\newcommand{\varint}{x}
\newcommand{\vareal}{y}
\newcommand{\basis}{B}
\newcommand{\Htmp}{H_{\textsc{tmp}}}
\newcommand{\diam}{d_{\mathcal{G}}}

\newcommand{\inneigh}{\mathcal{N}}
\newcommand{\Nag}{N}
\newcommand{\zetagen}{z}
\newcommand{\dimDZ}{d_Z}
\newcommand{\cvec}{c} 
\newcommand{\crho}{{e_1}}
\newcommand{\rhozi}{( \rho_I^{[i]}, z^{[i]} )}

\newcommand{\costval}{J}
\newcommand{\Npaths}{N_\gamma} %
\newcommand{\vectheta}{p_\theta} 
\newcommand{\vecthetaT}{p^\top_\theta} 
\newcommand{\vecnu}{q_\nu} 
\newcommand{\vecnuT}{q^\top_\nu} 
\newcommand{\vecgamma}{r_\gamma} 
\newcommand{\vecgammaT}{r^\top_\gamma} 

\newcommand{\LP}{\textsc{lp}}
\newcommand{\varzLP}{\varz^{\LP}} %
\newcommand{\varzLPell}{\varz^\LP_\indexa}%
\newcommand{\varintLP}{\varint^{\LP}}
\newcommand{\varintLPell}{\varint^\LP_\indexa} %
\newcommand{\varealLP}{\vareal^{\LP}}
\newcommand{\BLP}{B^{\LP}}

\newcommand{\indexa}{\ell}
\newcommand{\indexb}{m}
\newcommand{\indexc}{q}

\newcommand{\real}{\mathbb{R}}

\newcommand{\integer}{\mathbb{Z}}

\newcommand{\convS}{\text{conv}(P_I)}

\newcommand{\PP}{P}

\DeclarePairedDelimiter\ceil{\lceil}{\rceil}
\DeclarePairedDelimiter\floor{\lfloor}{\rfloor}

\newcommand{\subj}{\text{subj. to}}

\makeatletter
\newcommand{\StatexIndent}[1][3]{%
  \setlength\@tempdima{\algorithmicindent}%
  \Statex\hskip\dimexpr#1\@tempdima\relax}
\makeatother

\begin{document}

\title{Distributed Mixed-Integer Linear Programming\\ via Cut Generation
and Constraint Exchange}
    
\author{Andrea~Testa$^1$,
            Alessandro~Rucco$^1$, %
            and Giuseppe Notarstefano$^2$ %
\thanks{$^1$A. Testa and A. Rucco are with the
  Department of Engineering, Universit\`a del Salento, Lecce, Italy, 
  \texttt{name.lastname@unisalento.it.}  }
\thanks{$^2$G. Notarstefano is
    with the Department of Electrical, Electronic and Information Engineering, 
    University of Bologna, Bologna, Italy, \texttt{giuseppe.notarstefano@unibo.it.}}
\thanks{An early short version of this work appeared as \cite{testa2017dimilp}: the current article includes a 
much improved comprehensive treatment, an extended algorithm for general mixed-integer linear programs, 
new numerical computations, and an application to the multi-agent multi-task assignment problem.}
\thanks{
This result is part of a project that has received funding from
the European Research Council (ERC) under the European Union's
Horizon 2020 research and innovation programme (grant agreement
No 638992 - OPT4SMART).}}

\maketitle

\begin{abstract}
Many problems of interest for cyber-physical network systems can be
formulated as Mixed-Integer Linear Programs in which the constraints
are distributed among the agents. In this paper we propose a
distributed algorithmic framework to solve this class of optimization
problems in a peer-to-peer network with no coordinator and with limited
computation and communication capabilities. At each communication
round, agents locally solve a small linear program, generate suitable
cutting planes and communicate a fixed number of active constraints.
Within the distributed framework, we first propose an algorithm that,
under the assumption of integer-valued optimal cost, guarantees
finite-time convergence to an optimal solution. Second, we propose an
algorithm for general problems that provides a suboptimal solution up
to a given tolerance in a finite number of communication rounds. Both
algorithms work under asynchronous, directed, unreliable networks.
Finally, through numerical computations, we analyze the algorithm
scalability in terms of the network size. Moreover, for a multi-agent
multi-task assignment problem, we show, consistently with the theory,
its robustness to packet loss.
\end{abstract}

\IEEEpeerreviewmaketitle

\vspace{-0.02cm}

\section{Introduction}
\ac{MILPs} play an important role in several control problems, including control
of hybrid systems~\cite{bemporad1999control}, trajectory
planning~\cite{richards2002spacecraft}, and task
assignment~\cite{bellingham2003multi}.
For example, since non-convex functions can be approximated by means of
piecewise-linear functions, see, e.g.,~\cite{bemporad1999control}, 
nonlinear optimal control problems can be approximated by \ac{MILPs}.
Though \ac{MILPs} are known to be $\mathcal{NP}$-hard, numerous efficient
algorithms exist in a centralized setup. 
A widely used approach is the \emph{cutting-plane} method, see,
e.g.~\cite{gomory1963algorithm}, which is based on the iterative solution of
linear programming relaxations (obtained by neglecting integer constraints) 
and the generation of linear constraints (cuts) removing the current non-integer 
solution from the feasible set. 
In this paper, we consider the following \ac{MILP}
\begin{equation} \label{eq:MILP}
	\begin{split}
		\min_{\varz} &\; c^\top \varz\\ 
		\subj &\; a_i^\top \varz \leq b_i \,, \,i = 1, \ldots, n\\
		&\; \varz \in \integer^{d_Z} \times \real^{d_R}, \\
	\end{split}%
\end{equation}
where $d_Z$ and $d_R$ are the dimensions of the integer and real 
variables, $d=d_Z+d_R$, $a_i \in \real^{d}$, $b_i \in \real$, $c \in \real^{d}$, and $n$ is
the number of inequality constraints.

We deal with a distributed setup in which \ac{MILP} constraints are shared among
agents of a network, which aim at solving the whole problem by local
computations and communications. 
An important challenge that needs to be taken into account in a distributed
context is that the communication can be asynchronous, unreliable, and the
topology directed. Inspired by the centralized literature, we will propose a
distributed algorithmic framework based on cutting planes.

Only few works address the solution of MILPs in a purely distributed way, that is, with
peer processors communicating over a network, without the presence of a central (master) unit. 
For this reason we
organize the relevant literature to our paper in two main blocks: centralized
and parallel approaches to solve \ac{MILPs} in control applications, and
distributed algorithms solving Linear Programs (LPs) or convex programs arising
as relaxations or special versions of \ac{MILP}s.
As for centralized approaches, a Model Predictive Control (MPC) scheme to
solve constrained multivariable control problems is proposed in
\cite{bemporad2000piecewise,bemporad2002model}.
The MPC is formulated as a multi-parametric \ac{MILP} for which
solver~\cite{dua2000algorithm} is used.
In \cite{axehill2014parametric} a branch-and-bound procedure is devised for the
computation of optimal and suboptimal solutions to parametric \ac{MILP}s.
In~\cite{richards2002spacecraft} and~\cite{borrelli2006milp}, a collision-free
trajectory optimization problem for autonomous vehicles is formulated as a
\ac{MILP} solved by a branch-and-bound algorithm with branching
heuristics.
In \cite{chopra2015spatio}, a multi-robot routing problem
under connectivity constraints is shown to be formulated as an integer
program with binary variables, and then its LP relaxation is solved.
Recently, in \cite{takapoui2017simple} a heuristic based on the Alternating
Direction Method of Multipliers is used to approximately solve mixed-integer
linear and quadratic programs. The heuristic is applied to the
control of hybrid vehicles.
As for parallel methods, in \cite{kim2013scalable} a Lagrange relaxation
approach is used in order to decompose the overall \ac{MILP} into multiple
subproblems each of which is solved in a client-server parallel
architecture. The proposed solution is applied to the demand response control
in smart grids. 
In \cite{vujanic2016decomposition} a parallel dual decomposition method, relying
on a suitable tightening of the constraints, is proposed to approximately solve
structured \ac{MILPs} with local and coupling constraints.  
The algorithm is improved in \cite{falsone2017decentralized} by means of an
iterative tightening procedure. The methods are applied to charging control of
electric vehicles.

As for distributed optimization algorithms, we concentrate our review on schemes
solving linear programs or convex programs that represent a relaxation of common
mixed-integer programs.
In \cite{richert2016distributed_a} a robust, distributed algorithm is designed
to solve linear programs over networks with event-triggered communication.
In \cite{richert2016distributed_b} a distributed algorithm is proposed to find
valid solutions for the so called bargaining problem, which is an integer
program, by means of a linear program relaxation.
In \cite{fischione2011utility} a (mixed-integer) utility maximization problem is
addressed. The proposed solution is based on a convex relaxation obtained by
neglecting the integer constraint on the rates.
In \cite{wei2013distributed_alg} and \cite{wei2013distributed_conv} the authors
propose a Newton-type fast converging algorithm to solve Network Utility
Maximization problems with self-concordant utility functions.
In \cite{notarstefano2011distributed} constraints consensus algorithms are
proposed to solve abstract optimization programs (i.e., a generalization of LPs)
in asynchronous networks, while a distributed simplex algorithm is proposed in
\cite{burger2012distributed} to solve degenerate LPs and multi-agent assignment
problems.
A distributed version of the Hungarian method is proposed
in~\cite{chopra2017distributed} to solve LPs arising in multi-robot assignment problems.
In \cite{karaman2008large} and \cite{choi2009consensus} approximate
  solutions for task assignment (MILP) problems are proposed based respectively
  on a simplex ascent and an auction approach.
To conclude, \cite{franceschelli2013gossip,kuwata2011cooperative,falsone2018distributed} are
first attempts of proposing a distributed solution for \ac{MILP}s.
We discuss the main differences with our approach after the contributions paragraph.
The contributions of this paper are as follows. We propose a distributed
algorithmic framework, named Distributed Cutting Plane and Constraint Exchange
for \ac{MILP}s (\DiMILP/), based on the local generation of cutting planes, the
solution of local LP relaxations of \eqref{eq:MILP}, and the exchange of active
constraints.
Specifically, we propose two distributed algorithms. 
The first one, called \texttt{INT}-\DiMILP/, has guaranteed finite-time
convergence to an optimal solution of \eqref{eq:MILP} under the assumption of
integer-valued optimal cost. All relevant Integer Programs can, e.g., be casted
in this setup.
To remove the assumption of integer-valued optimal cost, we then propose an
algorithm, called $\epsilon$-\DiMILP/,\! ``practically'' \!solving~\eqref{eq:MILP},
i.e., computing a feasible point with cost exceeding the optimal one no more
than $\epsilon$.
Both algorithms involve, as local computations at each node, only a LP solver
(processing a number of constraints depending only on the problem dimension and
the number of neighboring nodes) 
and a simple procedure for the generation of cutting
planes. Thus, the algorithms are scalable in terms of local memory, computation,
and communication. 
Moreover, under slightly stronger assumptions on the graph, we provide a halting
condition allowing agents to stop the algorithm in a purely distributed way.
To the best of our knowledge, these are the first distributed algorithms
solving \ac{MILP}s in a distributed context. Notably, the proposed
algorithms work under asynchronous, unreliable, and directed networks, so that
they are immediately implementable in concrete scenarios. 
We highlight some meaningful differences with respect to the literature discussed above. 
In \cite{vujanic2016decomposition} and \cite{falsone2017decentralized},
suboptimal algorithms with performance guarantees for \ac{MILP}s are given and a
central unit is required.  
Papers \cite{franceschelli2013gossip}
and~\cite{kuwata2011cooperative} propose strategies to find feasible (suboptimal)
points for problem \eqref{eq:MILP}. 
Moreover, in \cite{kuwata2011cooperative}, agents perform the local computation
in a sequential order, while in \cite{franceschelli2013gossip} a gossip protocol
is considered in which one node per time becomes active. 
In~\cite{falsone2018distributed} a distributed algorithm with performance
guarantees, based on dual decomposition and a 
time-varying restriction technique, is proposed.
The algorithm applies to a class of optimization problems in which agents aim at minimizing the sum of local linear cost functions, 
subject to local linear constraints, while all local variables are coupled by global constraints. 
Our algorithm is developed for a different class of MILPs with common cost and local constraints. 
The algorithm in~\cite{falsone2018distributed} converges to a suboptimal solution with a performance guarantee that depends on the problem data. 
Our approach, instead, finds a suboptimal solution with an arbitrary tolerance and works over a general, possibly unreliable, asynchronous network.
Finally, although related approaches have been proposed in
\cite{notarstefano2011distributed} and \cite{burger2014polyhedral} for problems
with continuous decision variables, novel tools are needed in this paper due to
the mixed-integer nature of the problem.

The paper is organized as follows. 
In Section~\ref{sec:MILP} we recall the centralized cutting-plane approach for
\ac{MILPs}.
In Section~\ref{sec:DiMILP} we introduce a distributed meta-algorithm together
with two specific algorithms, while in Section~\ref{sec:analysis} we analyze
their convergence.
Numerical computations are provided in Section~\ref{sec:NumComp} for randomly
generated \ac{MILPs} and for a multi-agent multi-task assignment setup.

\paragraph*{Notation}
Given the decision variable $\varz\in\integer^{d_Z} \times\real^{d_R}$ of
\eqref{eq:MILP}, we denote by $\varint\in\integer^{d_Z}$ the vector of variables
subject to integer constraints, 
and by $\vareal\in\real^{d_R}$ the vector of variables
not required to be integer.
We denote by $e_\indexa$ the $\indexa$-th vector of the canonical basis (e.g.,
$e_1 = [1\, 0\, \ldots\, 0]^\top$) of proper dimension.
Given a vector $v \in \real^d$, we denote by $v_\indexa$ the $\indexa$-th
component of $v$.
Given two vectors $v, w \in \real^d$, $v$ is lexicographically greater than $w$,
$v >_{lex} w$, if there exists 
$\indexa \in \{1,\ldots, d\}$ such that $v_\indexa > w_\indexa$ and $v_m = w_m$
for all $m<\indexa$.
Given an inequality $a^\top \varz \leq b$ for $\varz\in\real^{d} $, with
$a\in\real^d$ and $b\in\real$, we use the following simplified notation
$\{ a^\top z \leq b\} := \{ z \in \real^d : a^\top z \leq b\}$ for the related
half-space.
The polyhedron induced by the inequality constraints
$a_i^\top \varz \leq b_i \,, \,i = 1, \ldots, n$, is
$P:=\bigcap_{i=1}^n \{ a_i^\top \varz \leq b_i \}$.
Recall that a polyhedron is a set described by the intersection of a finite
number of half-spaces.
In this paper we assume a LP solver is available.  In particular, we use the
simplex algorithm proposed in \cite{jones2007lexicographic} to find the unique 
\emph{lexicographically minimal optimal} (\emph{lex-optimal} for short) solution 
of degenerate linear programs.  
From now on, we call such a solver \textsc{LPlexsolv} and say that
it returns the lex-optimal solution of the solved LP. 
\textsc{LPlexsolv} also returns an optimal \emph{basis} identifying the lex-optimal solution.
Given a LP with constraint set $P:=\bigcap_{i=1}^n P_i$, with each $P_i$ a
half-space, a basis $B$ is the intersection of a \emph{minimal} number of
half-spaces $P_{\indexa_1}, \ldots,P_{\indexa_\indexc}$, $\indexc \leq d$, such
that the solution of the LP over the constraint set $B$ is the same as the one
over $P$. If the lex-optimal solution is considered, it turns out that $B$ is
the intersection of exactly $d$ half-spaces.

\section{Cutting-Plane framework for  Mixed-Integer Linear Programming}
\label{sec:MILP}
In this section we provide a brief description of one of the most used
centralized methods to solve a MILP, i.e., the cutting-plane approach, see
e.g.,~\cite{kelley1960cutting}. 
We introduce suitable cutting planes, namely intersection cuts and cost-based
cuts, and describe centralized algorithms for MILPs, based on these cuts, that
are relevant for our distributed framework.

\subsection{Cutting-Plane approach for MILP}\label{sec:CuttingPlane}
\begin{figure}[htpb]%
   \begin{center}
\includegraphics[width=.7\columnwidth]{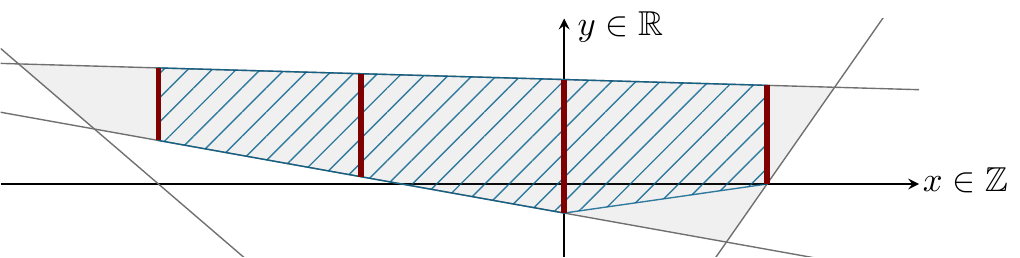}
   \caption{
   	Example of a polyhedron $P$ (gray area), with feasible set $P_I$ (solid red lines) with $\varint \in \integer$, $\vareal \in \real$ and its
   	convex hull (blue striped area).   
  }
	\label{fig:PSconvS}
   \end{center}
 \end{figure}

Let $P_I := P \cap (\integer^{d_Z} \times \real^{d_R})$ be the set of feasible
points of MILP~\eqref{eq:MILP}, also called \emph{mixed-integer set}. The
optimal solutions of \eqref{eq:MILP} are also optimal solutions of the following
LP, \cite{schrijver1998theory},
\begin{equation} \label{eq:MILP2}
\begin{split}
\min_{z} &\; c^\top \varz \\ 
\subj &\; \varz \in \convS \\
\end{split}
\end{equation}
where $\convS$ is the \emph{convex hull} of
$P_I$. A two-dimensional representation of a mixed-integer set and its convex hull is 
given in Figure~\ref{fig:PSconvS}.

It is worth noting that, if $P$ is a bounded polyhedron, by Meyer's
Theorem,~\cite{junger200950}, $\convS$ is a polyhedron.  
For this reason, we make the following assumption, which is
common in \ac{MILP} literature. 

\begin{assumption}[Boundedness and Feasibility] \label{ass:BoundFeas} %
  The polyhedron $P$ is bounded and $\convS$ is nonempty.~\oprocend
\end{assumption}

The main idea of the cutting-plane approach for \ac{MILP}s is to neglect the
  integer constraints on the decision vector $\varint$ (i.e.,
  $\varint \in \real^{d_Z}$), and iteratively solve \emph{relaxed} linear
  problems in which the polyhedron $P$ is tightened by additional half-spaces
  called \emph{cutting planes} or \emph{cuts}. The procedure terminates when the
  solution of the LP relaxation, call it
  $\varzLP = (\varintLP, \varealLP)$, is in $\convS$ with
  $\varintLP \in \integer^{d_Z}$.
  A valid cutting plane is a half-space containing $\convS$ but
  not $\varzLP$. 
We introduce a \textsc{CutOracle} subroutine defined as
follows. $\textsc{CutOracle}(\varzLP, P, c)$ returns the intersection of $p$
cutting planes, $\{\alpha^\top \varz \leq \beta \}$ with
$\alpha \in \real^{p\times d}$ and $\beta \in \real^p$, or $\real^d$.
An iterative cutting-plane scheme can be recast in the form of the following
meta-algorithm.

\bigskip

\textbf{Centralized Cutting-Plane Meta-Algorithm}: 
\begin{itemize}
\setlength\itemsep{0.2em}
\item[1.] \emph{Initialization}: set $P = \bigcap_{i=1}^n \{ a_i \varz \leq b_i \}$. 
\item[2.] \emph{LP solver}: Find an optimal solution
  $\varzLP = (\varintLP, \varealLP)$ of the LP relaxation of
  \eqref{eq:MILP} with polyhedron $P$.
\item[3.] \emph{Check feasibility}: if $\varintLP \in \integer^{d_Z}$, go to 6. 
\item[4.] \emph{Cutting-Plane}: $h = \textsc{CutOracle}(\varzLP, P, c)$.
\item[5.] \emph{Update}: $P = P \cap h$  and go to 2. 
\item[6.] \emph{Output}: $\varzLP$. 
\end{itemize}

\bigskip

It is worth noting that, at each iteration, the meta-algorithm uses
the entire set of inequality constraints, $P$, and all the cuts generated up to that iteration.
  
Following the meta-algorithm, numerous algorithms have been proposed in the
literature that have different convergence properties depending on \textsc{CutOracle} and 
the problem structure.
In the next subsections we describe two centralized algorithms. The key
ingredients of both algorithms are: i) \ac{MIG} cuts,
\cite{gomory1960algorithm}, and ii) \emph{cost-based cuts} generated according
to the current cost value $c^\top \varzLP$.

\subsection{Gomory's Cutting-Plane Algorithm}
Mixed-Integer Gomory (\ac{MIG}) cuts are cutting planes proposed in
\cite{gomory1960algorithm} for MILPs with integer-valued optimal cost.
Next, we introduce the notions of split disjunction, \cite{cook1990chvatal}, and
intersection cut, \cite{balas1971intersection}, that are intimately related to
\ac{MIG} cuts.

\begin{figure}[htpb]%
   \begin{center}
\includegraphics[scale=.7]{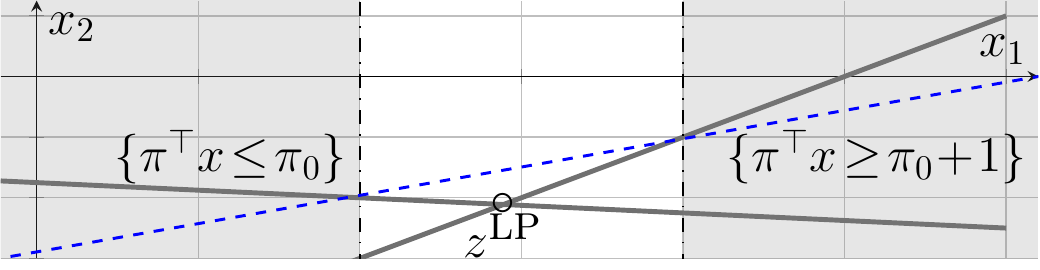}
	\caption{Example of split disjunction (gray shaded area) and intersection cut (dashed blue line) with respect to the basis $\BLP$ (solid lines) in $\real^2$. 
	}
	\label{fig:SplitIntersectionMIGcuts}
   \end{center}
\end{figure}

\begin{definition}[Split Disjunction \cite{cook1990chvatal}] \label{def:SplitDisj}
	Given $\pi \in \integer^{d_Z}$ and $\pi_0 \in \integer$, a split disjunction $D(\pi, \pi_0)$ is a set of the form 
	$D(\pi, \pi_0) := \{\pi^\top x\leq \pi_0\} \cup \{ \pi^\top \varint \geq \pi_0+1 \}$. \oprocend
\end{definition}
Let $\BLP$ be a basis for the optimal solution $\varzLP=(\varintLP,\varealLP)$ of 
a given LP relaxation of \eqref{eq:MILP}, and $D(\pi,\pi_0)$ a disjunction with respect 
to $\varintLP$. 
Let $C(\varzLP)$ be the translated simplicial cone\footnote{Given a cone
  $S\subset\real^d$ and a point $p\in\real^d$, the set $p + S$ is a translated
  cone with apex in $p$.} formed by the intersection of the half-planes defining
$\BLP$ with apex in $\varzLP$.
\emph{Intersection cuts} can be derived by considering the intersection between
the extreme rays of $C(\varzLP)$ and the hyperplanes defining the split 
disjunction $D(\pi,\pi_0)$. A more detailed definition can be found in
\cite{balas1971intersection}. %
A two dimensional representation of a split disjunction with  $\pi = [1 \; 0]^T$ 
is given in Figure~\ref{fig:SplitIntersectionMIGcuts} together with the
corresponding intersection cut 
(dashed blue line) 
for a given basis of $\varzLP$ (solid lines).

It can be shown that the \ac{MIG} cut with respect to $\varintLPell$, i.e., the
$\indexa$-th component of $\varintLP$, is the intersection cut to the split
disjunction $D(e_{\indexa},\floor{\varintLPell})$ and the basis $\BLP$, and is a
valid cutting plane~\cite{jorg2008k}.
Next we show how MIG cuts can be generated.
Let $\varzLP=(\varintLP,\varealLP)$ be the lex-optimal solution of a generic LP relaxation 
	of problem~\eqref{eq:MILP}
	with $\varintLP_{\ell} \not\in \integer$ for some $\ell \in\{1,\ldots,d_Z\}$, and let
	$\BLP = \{A_B \varz \leq b_B\}$ be an associated basis.
	Consider the split disjunction $D(e_\indexa, \floor{\varintLP})$
        which satisfies $\varzLP\not\in D(e_\indexa, \floor{\varintLP})$.  
	Let $r^\indexb$ be the $\indexb$-th column of $-A_B^{-1}$ and define $\lambda_\indexb$ as $\frac{\floor{\varzLPell} - e_\indexa^\top \varzLP}{e_\indexa^\top r^\indexb}$ 
  if $e_\indexa^\top r^\indexb < 0$, $\frac{\floor{\varzLPell} - e_\indexa^\top \varzLP +1}{e_\indexa^\top r^\indexb}$ 
  if $e_\indexa^\top r^\indexb > 0$ and $\infty$ if $e_\indexa^\top r^\indexb = 0$.
As it will be explained in Section~\ref{sec:analysis}, $r^\indexb$ represents an extreme ray (originating at $\varzLP$) intersecting the disjunction, 
while $\lambda_\indexb$ is the displacement along the extreme ray of the intersection point.
Let $\Lambda\in\real^d$ be the vector with $\indexb$-th entry $\Lambda_\indexb=1/\lambda_\indexb$ ($\Lambda_\indexb=0$ if $\lambda_\indexb=\infty$). 
Then, the intersection cut to the split disjunction $D(e_\indexa, \floor{\varintLPell})$ and the basis $\BLP$ is 
$h_{\textsc{mig}}=\{(\Lambda^\top A_B)^\top \varz \leq \Lambda^\top b_B-1\}$.

Let us consider now the first non-integer component of $\varintLP$, namely
$\varintLP_{k_{lex}}$ where $k_{lex} = \min \{ k = 1, \ldots, d_{Z} : \varintLP_k \notin
\integer \}$. 
We call \textsc{MIGoracle} the oracle that generates the \ac{MIG} cut $h_{\textsc{mig}}$  %
with respect to $\varintLP_{k_{lex}}$.
We highlight that \textsc{MIGoracle} returns $h_{\textsc{mig}}$ only when a variable subject to integer constraint is not integer, i.e., $\varintLP \notin \integer^{d_z}$.  
If the optimal solution $\varzLP$ obtained by solving the LP is such that $\varintLP \in \integer^{d_z}$, then \textsc{MIGoracle} does not return any constraint.

It is worth noting that if one implements the Centralized Cutting-Plane 
Meta-Algorithm using \textsc{MIGoracle} as \textsc{CutOracle}, the algorithm
may not converge. Indeed, a ``tailing-off'' phenomenon can be observed: a large number of cuts may be added
 without significant improvement in the cost.
A simple two-dimensional example is discussed in~\cite{owen2001disjunctive}.

 We are now ready to describe Gomory's Cutting-Plane Algorithm for MILPs with
 integer-valued optimal cost~\cite{gomory1960algorithm}. We provide a reformulation,
 given in \cite{jorg2008k}, for MILPs in the form \eqref{eq:MILP}.
 Given a basis $\BLP$ identifying the optimal solution $\varzLP$ and the
 corresponding cost function value $c^\top \varzLP$, we define \textsc{CutOracle}
 as follows
\begin{equation*} \label{eq:CutOracle_gomory}
	(h_{\textsc{mig}}, h_c) = \textsc{CutOracle}(\varzLP, \BLP, c)
\end{equation*}
where $h_{\textsc{mig}}$ is the \ac{MIG} cut generated by \textsc{MIGoracle} and
$h_c=\{c^\top \varz \geq \lceil c^\top \varzLP \rceil \}$.
It can be shown that, if the optimal objective function value is integer,
Gomory's cutting plane algorithm converges to an optimal solution in a finite number of
iterations~\cite{gomory1960algorithm}.

\subsection{Algorithms for $\epsilon$-suboptimal solutions}\label{sec:Alg_subopt}
Another relevant piece of literature regards algorithms solving MILPs up to an
arbitrary tolerance, namely based on the following notion of suboptimal
solution.

\begin{definition}[$\epsilon$-suboptimal solution]
\label{def:epsilon_suboptimal}
Given a MILP as in~\eqref{eq:MILP}, we say that $\varz^\epsilon$ is an
$\epsilon$-suboptimal solution to~\eqref{eq:MILP} if $\varz^\epsilon$ is
feasible and satisfies $c^\top \varz^\epsilon - c^\top \varz^\star \leq \epsilon$, where
$\varz^\star$ is an optimal solution to~\eqref{eq:MILP}. \oprocend
\end{definition}

In \cite{owen2001disjunctive} the authors propose a cutting-plane algorithm,
based on ``variable'' disjunctions, converging to an $\epsilon$-suboptimal
solution of the MILP in a finite number of iterations. In contrast to the usual
cutting-plane approaches, which generate constraints only at the optimal
solution of the current LP relaxation, the proposed algorithm generates
constraints at multiple, near-optimal vertices.
In \cite{jorg2008k} the author proposes an approximation algorithm which is not
a classical cutting-plane method because it generates cuts that might not be
valid for the mixed-integer set. Moreover, the algorithm relies on an inner
procedure to check if feasible points have been cut off, which involves the
solution of a MILP with integer-valued optimal cost. %

\section{A Distributed Cutting-Plane and Constraint Exchange Approach for
  MILPs} %
\label{sec:DiMILP}
Inspired by the centralized cutting-plane meta-algorithm, in this section we
first propose a distributed meta-algorithm, called \DiMILP/, based on the local
generation of cutting planes and the exchange of active
constraints.
As in the centralized case, it is not reasonable to provide a general, unified
convergence analysis for the meta-algorithm, but proper tools are needed for
specific algorithms. 
Thus, based on the high-level methodological approach, we provide an algorithm
for \ac{MILP}s with integer-valued optimal cost (as in
\cite{gomory1960algorithm}) and an approximation-based algorithm for general
problems.
Then, we discuss some key features of the proposed algorithms and provide a
distributed stopping criterion. 
We first formalize the distributed computation setup.

\subsection{Distributed Optimization Setup} 
In our distributed setup, we consider a network composed by a set of agents $V=\{1, \ldots, \Nag\}$. 
In general, the $n\geq \Nag$ constraints in problem~\eqref{eq:MILP} are
distributed among the agents, so that each agent knows only a small number of
constraints.
For simplicity, we assume one constraint $\{a_i^\top \varz \leq b_i\}$ is assigned
to the $i$-th agent, so that $\Nag = n$, but we will keep the two notations
separated to show that the algorithm can be easily implemented also when
$n>\Nag$ (i.e., more than one constraint is assigned to agents).  
The communication among the agents is modeled by a time-varying digraph
$\mathcal{G}_c(t)=(V,E(t))$, with $t\in\mathbb{N}$ being a universal slotted
time. A digraph $\mathcal{G}_c(t)$ models the communication in the sense that
there is an edge $(i,j) \in E(t)$ if and only if agent $i$ is able to send
information to agent $j$ at time $t$.  For each node $i$, the set of
\emph{in-neighbors} of $i$ at time $t$ is denoted by $\inneigh_i(t)$ and is the
set of all $j$ such that there exists an edge $(j,i) \in E(t)$.
A static digraph is \emph{strongly connected} if there exist a
directed path for each pair of agents $i$ and $j$.  
For a static digraph $\mathcal{G}_c=(V,E)$, we use $\diam$ to denote the graph diameter,
that is the maximum distance taken over all the pairs of agents $(i,j)$, where
the distance is defined as the length of the shortest directed path from $i$
to $j$. 
A time-varying digraph is \emph{jointly strongly connected} if, for all
$t \in \mathbb{N}$, $\cup_{\tau=t}^{\infty} \mathcal{G}_c(\tau)$ is strongly
connected.
The time-varying digraph is said to be \emph{(uniformly jointly) $L$-strongly
  connected} if there exists an integer $L \geq 1$ such that, for all
$t \in \mathbb{N}$, the graph $\cup_{\tau=t}^{t+L-1} \mathcal{G}_c(\tau)$ is
strongly connected.

\subsection{Meta-Algorithm description}
We now describe the proposed distributed meta-algorithm.
In contrast to the centralized approach, in the distributed setup, some agents
may need to solve local LP relaxations which are unbounded,
especially at the first iterations.
For this reason, we initialize the algorithm by assigning to each agent a set of
artificial constraints which are inactive for problem~\eqref{eq:MILP}. This
method is often referred to as big-M method.
Specifically, the decision variable of each agent is delimited by a box
constraint.  In particular, for a given, sufficiently
large $M>0$, we define the bounding box
$H_M:= \bigcap\limits_{\indexa=1}^{d} \left( \left\{ \varz_\indexa \leq{M} \right\} \cap
  \left\{ \varz_\indexa \geq{-M} \right\} \right) $. 

In the analysis we will need the following assumption. 
\begin{assumption}[Bounding Box] \label{ass:HM} 
  Given a MILP with polyhedron $\PP$ and a bounding box $H_M$, then
  $\PP \subseteq H_M$.~\oprocend
\end{assumption}

Each agent $i$ stores a fixed constraint $h^{[i]} = h_{i0} \cap H_M$, 
where each $h_{i0}$ is a local polyhedron (e.g., a single half-plane) known only 
by agent $i$,  
and updates two local states, namely $\zetagen^{[i]}$, associated to the decision
variable $\zetagen$ of the MILP, and $B^{[i]}$ being a candidate basis of the
problem.
At the generic (universal) time instant $t$, agent $i$ calls \textsc{CutOracle}
which returns, based on the current state $\zetagen^{[i]}(t)$, the \ac{MIG} cut,
$h_{\textsc{mig}}(t)$, obtained by \textsc{MIGoracle}, and a cost-based cut,
$h_{c}(t) = \{c^\top \varz \geq \sigma^{[i]}(t)\}$, where $\sigma^{[i]}(t)$ must
satisfy $\sigma^{[i]}(t) \geq c^\top \varz^{[i]}(t)$. 
Then, agent $i$ solves a local LP in which the common objective function
$\cvec^{\top}\zetagen$ is minimized subject to the following constraints: the
intersection of its neighbors' candidate bases,
$\bigcap_{j \in \inneigh_i(t)} B^{[j]}(t)$, its own candidate basis, $B^{[i]}(t)$,
the inequality constraint $h^{[i]}$, \ac{MIG} cut $h_{\textsc{mig}}(t)$, and the
cost-based cut $h_c(t)$.
This procedure is formalized in the following table. 

\begin{algorithm}[!ht]
\setstretch{1.3}
\renewcommand{\thealgorithm}{}
\floatname{algorithm}{Distributed Meta-Algorithm}
  \begin{algorithmic}
  \caption{\DiMILP/}\label{alg:distributed_meta}
    \State \textbf{State} $( \zetagen^{[i]}, B^{[i]} )$
    \State \textbf{Initialization}
    \StatexIndent[0.25]  $h^{[i]} = h_{i0} \cap  H_M $
    \StatexIndent[0.25]  $( \zetagen^{[i]}, B^{[i]}) =$ \Call{LPlexsolv}{$h^{[i]}$, {$\cvec$}}
    \State  \textbf{Evolution} 
    \StatexIndent[0.25]  $(h_{\textsc{mig}}(t), h_{\textsc{c}}(t)) =$ \Call{CutOracle}{$\zetagen^{[i]}(t)$, $B^{[i]}(t)$, $\cvec$}
    \StatexIndent[0.25]  $\Htmp(t)\!=\! \left( \bigcap_{j \in \inneigh_i(t)} \!B^{[j]}(t) \!\right) \cap B^{[i]}(t) \cap h^{[i]} \cap h_{\textsc{mig}}(t) \cap h_c(t)$
    \StatexIndent[0.25]  $(\zetagen^{[i]}(t+1), B^{[i]}(t+1)) =$ \Call{LPlexsolv}{$\Htmp(t)$, $\cvec$}
  \end{algorithmic}
\end{algorithm}

 Summarizing, at each communication round $t$, each agent sends to its out-neighbors a candidate basis, consisting of
$d$ linear constraints, and receives bases from its in-neighbors.
Then, each agent solves a local LP whose
number of constraints depends on the dimension $d$ and on the number of in-neighbors.
Notice that the distributed meta-algorithm does not require time synchronization: 
each agent can run its local routine (i.e., generate cuts and solve the local LP) at
its own rate, by directly using the available in-neighbor bases.

Next we provide two distributed algorithms. The first one provides an exact
solution under the assumption of integer-valued optimal cost. 
The second algorithm exploits a suitable reformulation and approximation of the
centralized problem to compute an $\epsilon$-suboptimal solution. 

\subsection{\texttt{INT}-\DiMILP/} \label{sec:ExactAlg}
We introduce a distributed algorithm to solve a MILP in the
form~\eqref{eq:MILP} under the assumption that the optimal cost is integer.
The algorithm was originally introduced in the preliminary conference
paper~\cite{testa2017dimilp}. 
Following the proposed distributed meta-algorithm, we set $h_{i0} = \{ a_i^\top \varz \leq b_i \}$. 
Then \textsc{CutOracle} consists of an oracle that generates \ac{MIG} cuts,
i.e., $h_{\textsc{mig}}(t) = \textsc{MIGoracle}(z^{[i]}(t), B^{[i]}(t))$, and another
oracle generating cost-based cuts, $h_c(t)=\{c^\top \varz \geq  \sigma^{[i]}(t) \}$, based on a 
simple ceiling of the current cost, that is, $\sigma^{[i]}(t) = \lceil c^\top \varz^{[i]}(t) \rceil$. 
We now state the convergence result for \texttt{INT}-\DiMILP/, whose proof is
given in the convergence analysis section. 
\begin{theorem}[\texttt{INT}-\DiMILP/ convergence]\label{thm:int_conv} %
  Let MILP~\eqref{eq:MILP} satisfy Assumptions~\ref{ass:BoundFeas} and~\ref{ass:HM}, 
  and let the optimal cost be integer-valued. 
  Assume agents communicate according to a jointly strongly connected
  communication graph, $\mathcal{G}_c(t)$, $t\geq0$, and run
  \texttt{INT}-\DiMILP/ distributed algorithm.
  Let $\costval^{[i]}(t) = \cvec^\top \zetagen^{[i]}(t)$ be the cost associated
  to the local candidate lex-optimal solution $\zetagen^{[i]}(t)$ of agent
  $i\in\{1,\ldots,\Nag\}$ at time $t\geq 0$.  
  Then, in a finite number of communication rounds the sequences
  $\{\costval^{[i]}(t)\}_{t \geq 0}$ and $\{\varz^{[i]}(t)\}_{t \geq 0}$,
  $i \in \{ 1, \ldots \Nag\}$, converge, respectively, to the (integer-valued) optimal cost and to
  the lex-optimal solution of problem~\eqref{eq:MILP}.\oprocend
\end{theorem}

Consistently with the (centralized) Gomory's Cutting-Plane Algorithm, \texttt{INT}-\DiMILP/ 
requires an integer-valued optimal cost. 
Such assumption is needed in order to guarantee that cost-based cuts 
(obtained by rounding up the current local cost value $c^\top z^{[i]}(t)$) 
are valid cutting planes. 
One way to guarantee integer-valued optimal cost is to take into account \ac{MILP} instances 
in which the cost function value depends only on the integer variables and the cost vector 
has rational components (which, through a suitable scaling, is equivalent to assuming 
they are integer valued). 
We point out that for a variety of practical problems, as, for example, 
scheduling, cutting stock, warehouse location,~\cite{junger200950}, 
the optimal cost is integer because all the decision variables are required to be integer and 
the coefficients of the cost vector $c$ are rational.
On the other hand, in many cases of interest, an
integer-valued optimal cost cannot be guaranteed a-priori. In the next subsection we remove
this assumption and propose a distributed algorithm based on a suitable
reformulation and approximation of the centralized problem.

\subsection{$\epsilon$-\DiMILP/} \label{sec:ApproxAlg}
Here we propose a distributed algorithm, based on \DiMILP/ meta algorithm, which 
computes an $\epsilon$-suboptimal solution for general \ac{MILP}s. 
We first consider an equivalent formulation of MILP~\eqref{eq:MILP}, namely the
epigraph form %
\begin{equation} \label{eq:MILP_epigraph}
	\begin{split}
		\min_{\rho, \varz} &\; \rho\\ 
		\subj &\; a_i^\top \varz \leq b_i \,, \,i = 1, \ldots, n\\
		&\; c^\top \varz \leq \rho\\
		&\; \rho \in \real, \varz \in \integer^{d_Z} \times \real^{d_R} \\
	\end{split} 
\end{equation}
where $\rho \in \real$ is a new decision variable. 
Problems~\eqref{eq:MILP} and~\eqref{eq:MILP_epigraph} are equivalent in the
sense that $(\rho^\star, \varz^\star)$ is optimal for the epigraph
form~\eqref{eq:MILP_epigraph} if and only if $\varz^\star$ is optimal for
problem~\eqref{eq:MILP} and $\rho^\star = c^\top \varz^\star$.

Now, we make a change of variables $\epsilon \rho_I = \rho$, with $\epsilon>0$,
and approximate problem \eqref{eq:MILP_epigraph} by constraining the new
variable $\rho_I$ to be integer. 
The resulting approximate MILP is
\begin{equation} \label{eq:MILP_epigraphrhoI}%
	\begin{split}
		\min_{\rho_I, \varz} &\; \rho_I \\ 
		\subj &\; a_i^\top \varz \leq b_i \,, \,i = 1, \ldots, n\\
		&\; c^\top \varz \leq \epsilon \rho_I \\
		&\; \rho_I \in \integer, \varz \in \integer^{d_Z} \times \real^{d_R}, \\
	\end{split} %
\end{equation}
where the constant $\epsilon>0$ is neglected in the objective function since it
does not affect the minimization. %

With problem~\eqref{eq:MILP_epigraphrhoI} at hand, we are ready to introduce
$\epsilon$-\DiMILP/ distributed algorithm. %
The algorithm is obtained by properly targeting \DiMILP/ distributed
meta-algorithm to problem~\eqref{eq:MILP_epigraphrhoI}. That is, we set
$(\rho_I, \varz)$ as extended decision variable, $\crho$ as cost
(as opposed to $\cvec$), and set
$h_{i0} = \{ a_i \varz \leq b_i \} \cap \{ c^\top \varz \leq \epsilon \rho_I
\}$ as local constraint.
The polyhedron induced by the inequality constraints is
$P^\epsilon = \{ (\rho_I, \varz) \in \real^{d+1} : a_i^\top \varz \leq b_i, i
= 1,\ldots,n, \, c^\top \varz \leq \epsilon \rho_I\}$.
We denote $H_M^\epsilon$ the bounding box associated to the
decision variable $(\rho_I, z)$. 
Consistently with the distributed meta-algorithm, each agent generates a
\ac{MIG} cut
$h_{\textsc{mig}}(t) = \textsc{MIGoracle}((\rho_I^{[i]}(t), z^{[i]}(t)),
B^{[i]}(t))$
and a cost-based cut $h_c(t) = \{ \rho_I \geq \sigma^{[i]}(t)\}$ where
$\sigma^{[i]}(t) = \lceil \rho_I^{[i]}(t) \rceil$. Then, it solves a local LP
based on the generated cuts and on the neighboring bases.
To better highlight the connection with the meta-algorithm, a pseudo-code
description of $\epsilon$-\DiMILP/ is reported in the following table.

\begin{algorithm}[htpb]
\setstretch{1.3}
\renewcommand{\thealgorithm}{}
  \floatname{algorithm}{Distributed Algorithm}
\caption{$\epsilon$-\DiMILP/}\label{alg:distributed_epsilon}
  \begin{algorithmic}[0]		
    \State \textbf{State} $( \rhozi, B^{[i]} )$
    \State \textbf{Initialization}
    \StatexIndent[0.25]  $h^{[i]} = h_{i0} \cap H_M $
    \StatexIndent[0.25]  $( \rhozi, B^{[i]}) =$ \Call{LPlexsolv}{$h^{[i]}$, $\crho$}
    \State  \textbf{Evolution}
    \StatexIndent[0.25]  $h_{\textsc{mig}}(t) =$ \Call{MIGoracle}{$(\rho_I^{[i]}(t), z^{[i]}(t))$, $B^{[i]}(t)$}
    \StatexIndent[0.25]  $h_{\textsc{c}}(t) = \{ \rho_I \geq \lceil \rho^{[i]}_I(t) \rceil \}$
    \StatexIndent[0.25]  $\Htmp(t)\!=\! \left(\! \bigcap_{j \in \inneigh_i(t)} \!B^{[j]}(t)\! \right) \! \cap B^{[i]}(t) \cap h^{[i]} \cap h_{\textsc{mig}}(t) \cap  h_c(t)$
    \StatexIndent[0.25]  $( (\rho_I^{[i]}\!(t\!+\!1),\! z^{[i]}\!(t\!+\!1)),\! B^{[i]}\!(t\!+\!1) )\!\! =\!\!$ \Call{LPlexsolv}{$\Htmp(t)$,\! $\crho$}
  \end{algorithmic}
\end{algorithm}

The convergence properties of $\epsilon$-\DiMILP/ are stated in the next
theorem whose proof is given in the analysis section. 
\begin{theorem}[$\epsilon$-\DiMILP/ convergence]\label{thm:eps_conv} %
  Let MILP~\eqref{eq:MILP} satisfy Assumptions~\ref{ass:BoundFeas}
  and~\ref{ass:HM} with 
  $M \geq \max \left( -\min_{z\in P}c^\top z/ \epsilon, \ceil{ \max_{z\in
        P}c^\top z / \epsilon} \right)$.
  Assume agents communicate according to a jointly strongly connected
  communication graph, $\mathcal{G}_c(t)$, $t\geq0$, and run $\epsilon$-\DiMILP/
  distributed algorithm. 
  Then, in a finite number of communication rounds the sequences
  $\{\varz^{[i]}(t))\}_{t \geq 0}$, $i \in \{ 1, \ldots \Nag\}$, converge to an
  $\epsilon$-suboptimal solution of~\eqref{eq:MILP}.~\oprocend
\end{theorem}

From Theorem~\ref{thm:eps_conv} it follows immediately that if the optimal cost of MILP~\eqref{eq:MILP} consists of $q$ decimal
  digits, then, by setting $\epsilon = 10^{-q}$, 
  agents compute an optimal solution of~\eqref{eq:MILP} in a finite number of
  communication rounds.
  \begin{remark}(Multiple Cuts)\label{rmk:multi_cut}
  Both \texttt{INT}-\DiMILP/ and $\epsilon$-\DiMILP/ can be implemented by generating multiple \ac{MIG}
  cuts at the generic communication round $t$. That is, together with the MIG cut with respect to the first non-integer component
  generated by \textsc{MIGoracle}, agent $i$ also generates intersection cuts for other non-integer 
  entries. %
  The introduction of multiple cuts can be a useful tool for faster and
  numerically robust convergence, as shown in the numerical computations.
  Proofs are provided for a single cut, but the extension to multiple
  cuts follows by using similar arguments. \oprocend
  \end{remark}

\subsection{Discussion of Main Algorithm Features}%
  We discuss some interesting features of the proposed
  algorithms.
  First, our distributed algorithms only require the communication graph to be
  jointly strongly connected, and the universal time does not need to be know by
  the agents. Indeed, agents do not use it in the local updates. This implies,
  as we will show in the numerical computations, that the algorithms work under
  asynchronous and unreliable communication networks and, in particular, in
  networks subject to packet loss.
  Second, agents can realize that convergence has occurred, and thus halt the
  algorithm, in a purely distributed way under slightly stronger assumptions on
  the graph.
We recall that, for a static digraph, the diameter $\diam$ is the maximum
distance taken over all the pairs of agents ($i,j$), where the distance is
defined as the length of the shortest directed path from $i$ to $j$.
Since our distributed algorithms converge in a finite number of communication
rounds, it can be shown that, for static communication digraphs, each agent
running the algorithm can stop the algorithm 
if its basis has not changed for $2\diam +1$ communication rounds, see,
e.g.,~\cite{notarstefano2011distributed}. 
It is worth noting that, in the initialization step, each agent can compute the
graph diameter by a simple flooding algorithm. 
By similar arguments it can be shown that, if the graph is (uniformly jointly)
$L$-strongly connected, then each agent can stop the algorithm if the value of
its basis has not changed after $2 L \Nag+1$ communication rounds. 
Third, the distributed meta-algorithm (so as the two specific algorithms)
involves local computations and communications that depend on the dimension $d$ of the decision
variable and on the number of in-neighbors.
  Indeed, an agent sends to neighbors a candidate basis, which is a collection
  of $d$ linear constraints, generates cutting planes based on simple local
  computations, and solves a LP.
  Thus, the main computational burden for the $i$-th agent is due to the
  solution of a LP (solvable in polynomial time) with $d$ variables and
  $d\times \inneigh_i(t)$ constraints.

  Finally, we conclude this discussion by highlighting that the idea of solving
  a suitable epigraph approximation of the original problem, proposed for
  $\epsilon$-\DiMILP/, can be used also in a centralized setup. This would
  result into a pure cutting plane algorithm providing an $\epsilon$-suboptimal
  solution for general \ac{MILP}s with a computation burden comparable with the
  Gomory algorithm (which however works only under the assumption of
  integer-valued optimal cost).
  On this regard, we point out that other approximate algorithms, as the ones
  mentioned in Section~\ref{sec:Alg_subopt}, require inner, computationally
  expensive, procedures to find a suboptimal solution.

\section{Convergence Analysis of  \DiMILP/ algorithms}
\label{sec:analysis}
In this section we provide a twofold result. First, we provide two technical
lemmas that hold for the \DiMILP/ distributed meta-algorithm and that can be
used to prove the convergence of a class of distributed optimization algorithms
based on this approach.
Specifically, we prove that for each agent the local cost and
the local state converge in a finite number of communication rounds if the
generation of cost-based cuts is assumed to stop in finite time. Then, we prove
that under the same condition, consensus among all the agents is attained for
the costs and for the candidate lex-optimal solutions.
Second, relying on these results, we prove the convergence of
\texttt{INT}-\DiMILP/ and $\epsilon$-\DiMILP/ distributed algorithms. Namely,
for \texttt{INT}-\DiMILP/ we show that agents agree on the lex-optimal solution
of \ac{MILP}~\eqref{eq:MILP} under the assumption of integer-valued optimal
cost. For $\epsilon$-\DiMILP/ we show that they agree on an
$\epsilon$-suboptimal solution of \ac{MILP}~\eqref{eq:MILP}.

\subsection{Property of \ac{MIG} Cuts to a Basis}
Here, we prove a property, that will be used in our convergence analysis,
holding for \ac{MIG} cuts generated with respect to a basis (rather than with
respect to the entire polyhedron as usually done in centralized algorithms).
\begin{lemma}\label{prop:BasisInterCuts}
	Let $\varzLP$ be the lex-optimal solution of a generic LP relaxation 
	of problem~\eqref{eq:MILP} 
	with $\varzLPell \not\in \integer$ for some $\ell \in\{1,\ldots,d_Z\}$, and 
	$\BLP = \{A_B \varz \leq b_B\}$ an associated basis. 
	Let $D(e_\indexa, \floor{\varzLPell})$ be the split disjunction
        $\{e_\indexa^\top \varint \leq \floor{\varzLPell} \} \cup
        \{e_\indexa^\top \varint \geq \floor{\varzLPell}+1\}$,
        which does not contain $\varzLP$.  
	Let $h:= \{ \alpha^\top \varz \leq \beta \}$ be the intersection cut to
        the disjunction $D(e_\indexa, \floor{\varzLPell})$ and to the basis
        $\BLP$. Then, the lex-optimal solution, $\varz^B$, obtained by
        minimizing the linear cost $c^\top \varz$ over $\BLP \cap h$ is such
        that its $\indexa$-th component is either
        $\varz^B_\indexa = \floor{\varzLPell}$ or
        $\varz^B_\indexa = \floor{\varzLPell} + 1$.
\end{lemma}
\begin{IEEEproof}
  The intersection of the half-planes defining $\BLP$ defines a translated cone
  $C(\varzLP)$ with apex $\varzLP$. %
  Points along an extreme ray\footnote{A nonzero vector $r$ of a polyhedral cone
    $C$ is called an extreme ray if there are $d - 1$ linearly independent
    constraints that are active at $r$, \cite{bertsimas1997introduction}.}
  $r^\indexb \in \real^d$, $\indexb \in \{1, \ldots, d\}$, associated to
  $C(\varzLP)$ %
  are described by $\varz = \varzLP + \mu r^\indexb$, $\mu \in \real$ and
  $\mu \geq 0$.
  Thus, vectors $\varz \in \BLP$ can be described by a positive linear
  combination of the extreme rays $r^\indexb$ of the basis $\BLP$, i.e.,
  $\varz = \varzLP + \sum_{\indexb = 1}^d \mu_\indexb r^\indexb, \mu_\indexb \geq 0$. 
  For those $r^\indexb$ such that $e_\indexa^\top r^\indexb \neq 0$, let us
  define $\varz^\indexb = \varzLP + \lambda_\indexb r^\indexb$ with
\begin{equation*}\label{eq:intercut_lambda}
	\lambda_\indexb = 
	\begin{cases}
	\frac{\floor{\varzLPell} - e_\indexa^\top \varzLP}{e_\indexa^\top r^\indexb} & \text{if} \,\,\, e_\indexa^\top r^\indexb < 0\\
	\frac{\floor{\varzLPell} - e_\indexa^\top \varzLP +1}{e_\indexa^\top r^\indexb} & \text{if}\,\,\, e_\indexa^\top r^\indexb > 0.\\
	\end{cases}
\end{equation*}
Then $\varz^\indexb \in \BLP$ because $\lambda_\indexb \geq 0$. Moreover, 
$\varz^\indexb$ is the point obtained as intersection of the ray $r^\indexb$
with the disjunctive hyperplanes $e_\indexa^\top \varz = \floor{\varzLPell}$ or
$e_\indexa^\top \varz = \floor{\varzLPell}+1$.  In fact, by solving
\[
\begin{cases}
	& \varz = \varzLP + \mu r^\indexb \\
	& e_\indexa^\top \varz = \floor{\varzLPell} +1 \\
\end{cases} \,,
\]
we have
\[
\begin{split}
& e_\indexa^\top(\varzLP+\mu r^\indexb) = \floor{\varzLPell}+1 \\
& \mu e_\indexa^\top r^\indexb = \floor{\varzLPell}+1 - e_\indexa^\top \varzLP\\ %
& \mu = \frac{\floor{\varzLPell}+1 - e_\indexa^\top \varzLP}{e_\indexa^\top
  r^\indexb}, 
\end{split}
\]
which is well defined and nonnegative if $e_\indexa^\top r^\indexb>0$. A similar
argument holds for $e_\indexa^\top \varz = \floor{\varzLPell}$ which gives the
condition $e_\indexa^\top r^\indexb<0$. 
If $e_\indexa^\top r^\indexb = 0$, the extreme ray $r^\indexb$ is parallel to
the disjunctive hyperplanes $e_\indexa^\top \varz = \floor{\varzLPell}$ and
$e_\indexa^\top \varz = \floor{\varzLPell}+1$, \cite{andersen2005split}. 
The intersection cut to the split disjunction $D(e_\indexa, \floor{ \varzLPell})$ and the basis $\BLP$, 
$h = \{ \alpha^\top \varz \leq \beta \}$, is then defined by the hyperplane passing through the intersection 
points $\varz^\indexb$ such that $e_\indexa^\top r^\indexb \neq 0$ and does not intersect the rays such that $e_\indexa^\top r^\indexb = 0$.  
Now, the lexicographic minimization of the linear cost function over $h \cap \BLP$ returns $\varz^B$ 
(cutting off $\varzLP$) and the new basis characterized by $d-1$ constraints from $\BLP$ and $h$. 
The intersection of the half-planes defining the new basis defines a translated cone $C(\varz^B)$ with apex $\varz^B$. 
Let us consider the extreme ray $r^{\bar{\indexb}}$ defined by the $d-1$
constraints of $\BLP$ in the new basis. Now, $r^{\bar{\indexb}}$ cannot be
parallel to the intersection cut $h$ otherwise the associated $d-1$ constraints
could not belong to the new basis. Thus, there exists a point
$\varz^{\bar{\indexb}}$ on the ray intersecting one of the two disjunctive
hyperplanes.
By construction, $\varz^{\bar{\indexb}}$ belongs also to the intersection cut
$h$ and thus it must be the unique intersection point $\varz^B$.
Therefore, $\varz^B$ has $\ell$-th component equal to $\floor{\varzLPell}$ or $\floor{\varzLPell}+1$.
\end{IEEEproof}

\subsection{Technical Results for the \DiMILP/ Meta-Algorithm}

\begin{lemma}[Local convergence] \label{lemma:conv} 
  Let MILP~\eqref{eq:MILP} satisfy Assumptions~\ref{ass:BoundFeas}
  and~\ref{ass:HM}.  Assume agents run an algorithm based on \DiMILP/ distributed
  meta-algorithm such that, at each node, \textsc{CutOracle} generates
  cost-based cuts $\{ c^\top\varz \geq \sigma^{[i]}(t)\}$ with $\sigma^{[i]}(t)$
  taking values in a finite set.
  Then, in a finite number of
  communication rounds, for all $i\in\{1,\ldots,\Nag\}$,
\begin{enumerate}
  \item[i)] the sequence $\{\costval^{[i]}(t)\}_{t \geq 0}$ converges to a constant
  value $\bar{\costval}^{[i]}$, and
  \item[ii)] the sequence $\{\zetagen^{[i]}(t)\}_{t \geq 0}$ converges to a feasible 
  $\bar{\zetagen}^{[i]}$ (with $\bar{\zetagen}^{[i]}_1, \ldots,
  \bar{\zetagen}^{[i]}_{\dimDZ}$ integer). 
\end{enumerate}
\end{lemma}
 \begin{IEEEproof}
To prove that sequence $\{\costval^{[i]}(t)\}_{t \geq 0}$ converges
to a constant $\bar{\costval}^{[i]}$ in a finite number of communication rounds we proceed in four steps. 

First, the sequence $\{\costval^{[i]}(t)\}_{t \geq 0}$ is monotonically
non-decreasing. Indeed, to compute $\basis^{[i]}(t+1)$, the $i$-th agent
minimizes the common objective function subject to the constraint set $\Htmp(t)$
which, by construction, is a subset of the basis $\basis^{[i]}(t)$. 
Second, the sequence is bounded from above by the optimal cost
${\costval}^\star$ to~\eqref{eq:MILP} (since $\Htmp(t)$ always includes
$\convS$). 
Thus, $\costval^{[i]}(t)$ converges to some
$\bar{\costval}^{[i]}\leq {\costval}^\star$.
To conclude the proof of statement i), we recall that at each $t$,
$h_{c}(t) = \{c^\top \varz \geq \sigma^{[i]}(t)\}$, with
$\sigma^{[i]}(t) \geq J^{[i]}(t)$ and by assumption $\sigma^{[i]}(t)$ takes
value in a finite set.

To prove the second statement, first notice that we have just proved that there
exists a time $t_0$ such that
$\costval^{[i]}(t) = \bar{\costval}^{[i]} \,, \forall t \geq t_0$.
Now let us consider the sequence of the first component of the state associated to the integer
decision variable, i.e., $\{\zetagen_1^{[i]}(t)\}_{t \geq t_0}$. %
The sequence $\{\zetagen_1^{[i]}(t)\}_{t \geq t_0}$ is non-decreasing because the local cost 
value is constant after $t_0$, $\Htmp(t)$ is a subset of $B^{[i]}(t)$,
$\forall t \geq t_0$, and the sequence is constructed by taking into account the
lex-optimal solution of the local problem.
Moreover, it is upper bounded by $M$, and
therewith convergent with limit $\tilde{\zetagen}^{[i]}_1$. %
So, there exists a time $t_1 \geq t_0$ such that 
$\lceil \tilde{\zetagen}^{[i]}_1 \rceil - 1 < \zetagen^{[i]}_1(t) \leq \lceil \tilde{\zetagen}^{[i]}_1\rceil $, $\forall t \geq t_1$. 
Following the evolution of the meta-algorithm, to compute
$\zetagen_1^{[i]}(t_1+1)$, first agent $i$ generates, through
\textsc{CutOracle}, a \ac{MIG} cut $h_{\textsc{mig}}(t_1)$ to the split
disjunction $D(e_1, \lceil{ \tilde{\zetagen}_1^{[i]} } \rceil -1)$ and to the
current basis $B^{[i]}(t_1)$.
Then it collects the constraints from its neighbors, builds up $\Htmp(t_1)$, and
calls \textsc{LPlexsolv}. This returns a new lex-optimal solution
$\zetagen_1^{[i]}(t_1+1)$, which is greater than or equal to the solution
obtained by minimizing over $B^{[i]}(t_1) \cap h_{\textsc{mig}}(t_1)$, since
$\Htmp(t_1)$ is a subset of it.
By Lemma~\ref{prop:BasisInterCuts}, minimizing the linear cost function
over the set of the current basis and the intersection cut is such that the
first component is either $\lceil{ \tilde{\zetagen}_1^{[i]} } \rceil -1$ or
$\lceil{ \tilde{\zetagen}_1^{[i]} } \rceil$.
But, being $\lceil \tilde{\zetagen}^{[i]}_1 \rceil - 1 < \zetagen^{[i]}_1(t)$ for all
$t\geq t_1$, it must hold $\zetagen_1^{[i]}(t_1+1) = \lceil{
  \tilde{\zetagen}_1^{[i]} } \rceil$ and, thus, $\zetagen_1^{[i]}(t) = \lceil{
  \tilde{\zetagen}_1^{[i]} } \rceil$ for all $t>t_1+1$.
Therefore, we have shown that the sequence $\{\zetagen^{[i]}_1(t)\}_{t \geq t_0}$
converges to $\bar{\zetagen}^{[i]}_1 = \lceil \tilde{\zetagen}^{[i]}_1 \rceil$ in a
finite number of communication rounds.
The same argument can now be applied to the remaining integer 
components of the decision vector, 
$\zetagen^{[i]}_2, \ldots, \zetagen^{[i]}_{\dimDZ}$. %
So there exists a time $t_{\dimDZ}$ such that 
$(\zetagen^{[i]}_1(t), \ldots, \zetagen^{[i]}_{\dimDZ}(t)) = (\bar{\zetagen}^{[i]}_1,
\ldots, \bar{\zetagen}^{[i]}_{\dimDZ})$, $\forall t \geq t_{\dimDZ}$, and thus
agent $i$ will not generate cutting planes anymore. 
Let $T_{\dimDZ}$ be such that $(\zetagen_1^{[i]}(t), \dots, \zetagen_{\dimDZ}^{[i]}(t)) \in \integer^{\dimDZ}$,
$\forall t \geq T_{\dimDZ}$, $\forall i \in \{1, \ldots, \Nag\}$. For
$t \geq T_{\dimDZ}$, no cutting planes will be generated in the network. 
This means that the number of possible different bases that $i$-th agent can
receive from its neighbors is finite (specifically, the combination of all the
constraints in the network at $T_{\dimDZ}$).
Therefore, since the lexicographic cost is nondecreasing and due to the finite
number of possible bases, also the non-integer variables
$(\zetagen_{\dimDZ+1}^{[i]}(t), \ldots, \zetagen_{d}^{[i]}(t))$ will
converge in a finite number of communication rounds, thus concluding the
proof.
\end{IEEEproof}

Next we prove that agents reach consensus on (a common) cost value and solution
estimate in a finite number of communication rounds.
The proof follows similar arguments as in \cite{notarstefano2011distributed} but
we report all the steps for the sake of
completeness.

\begin{lemma}[Consensus] \label{lemma:agreement} 
  Assume that the sequences $\{J^{[i]}(t)\}_{t\geq0}$ and
  $\{\zetagen^{[i]}(t)\}_{t \geq 0}$ defined as in Lemma~\ref{lemma:conv}
  converge to a constant value $\bar{J}^{[i]}$ and a feasible
  $\bar{\zetagen}^{[i]}$, respectively, $\forall i \in \{1, \ldots, \Nag\}$.
  Assume the communication network, $\mathcal{G}_c(t)$, is jointly strongly
  connected. Then, $\bar{\costval}^{[i]} = \bar{\costval}^{[j]}$ and
  $\bar{\zetagen}^{[i]}=\bar{\zetagen}^{[j]}$ for all
  $i,j \in \{ 1, \ldots, \Nag\}$.
\end{lemma}
 \begin{IEEEproof}
  We start by proving the consensus in cost. 
  Suppose, by contradiction, that two cost sequences, associated with two
  agents, say $i$ and $j$, converge to two different values
  $\bar{\costval}^{[i]}$ and $\bar{\costval}^{[j]}$, respectively. Let, without
  loss of generality, $\bar{\costval}^{[j]} > \bar{\costval}^{[i]}$.
  In Lemma~\ref{lemma:conv} we
  have shown that both sequences, $\{ \costval^{[i]}(t) \}_{t \geq 0}$ and
  $\{ \costval^{[j]}(t) \}_{t \geq 0}$, are monotonically non-decreasing and
  convergent in a finite number of communication rounds.  
  Since $\bar{\costval}^{[j]} - \bar{\costval}^{[i]} > 0$, 
  there must exist
  $T_{\delta_{ij}}>0$ such that
  $\bar{\costval}^{[j]} \geq \costval^{[j]}(t) >\bar{\costval}^{[i]} \geq \costval^{[i]}(t)$, $\forall t>T_{\delta_{ij}}$. 
  Now since the communication graph is jointly strongly connected,
  for each time $t\geq 0$ and each pair of agents $(i,j)$, there exists a
  sequence of nodes $\{l_1,\ldots,l_k\}$ and an increasing sequence of time
  instants $\{t_{l_1},\ldots,t_{l_{k}+1}\}$, with
  $t \leq t_{l_1} < \ldots < t_{l_{k}+1}$, such that the directed edges
  $\{(j,l_1),(l_1,l_2),\ldots,(l_k,i)\}$ belong to the digraph at times
  $\{ t_{l_1},\ldots,t_{l_{k}+1}\}$,
  \cite{notarstefano2011distributed}.
  To compute $\basis^{[l_1]}(t_{l_1}+1)$, 
  agent $l_1$ minimizes the common objective function subject to the constraint
  set $\Htmp(t_{l_1})$ which, by construction, is a subset of the basis
  $\basis^{[j]}(t_{l_1})$. %
  Therefore, $\costval^{[l_1]}(t_{l_1}+1) \geq \costval^{[j]}(t_{l_1})$.
  Iterating, we have $\costval^{[i]}(t_{l_{k}+1}+1) \geq \costval^{[j]}(t_{l_1})$.
  Therefore, for each $t>T_{\delta_{ij}}$ there exists $\tau>0$ such that
  $\costval^{[i]}(t+\tau) \geq \costval^{[j]}(t)$, 
  which leads to a contradiction, thus proving that
  $\bar{\costval}^{[1]} = \ldots= \bar{\costval}^{[\Nag]}$. %

  To prove that
  $\bar{\zetagen}^{[i]} = \bar{\zetagen}^{[j]} \, \forall \, i, j \in \{ 1,
  \ldots, \Nag\}$,
  first recall that each sequence $\{\zetagen^{[i]}(t)\}_{t \geq 0}$ converges to
  $\bar{\zetagen}^{[i]}$ $\forall i \in \{1, \ldots, \Nag \}$ in a finite number
  of communication rounds. 
Let us suppose, by contradiction, that agents do not reach consensus on the
decision variable, then there exist two agents, say $i$ and $j$, such that
$\bar{z}^i >_{lex} \bar{z}^j$.
Using again the assumption that the graph is jointly strongly connected, there
exists a directed path $\{(i,l_1),(l_1,l_2),\ldots,(l_k,j)\}$ at times
$\{ t_{l_1},\ldots,t_{l_{k}+1}\}$. Since consensus has been reached on the cost,
it must hold (by applying to the lexicographic ordering the same argument used
for the cost) that
$\bar{z}^j=\zetagen^{[j]}(t_{l_{k}+1}+1 ) \geq_{lex}
\zetagen^{[l_k]}(t_{l_{k}+1} ) \geq_{lex} \ldots \geq_{lex}
\zetagen^{[i]}(t)=\bar{z}^i$
which leads to a contradiction, thus concluding the proof.
\end{IEEEproof}	

\subsection{Proof of Theorem~\ref{thm:int_conv}}
  We apply Lemma~\ref{lemma:conv} to MILP~\eqref{eq:MILP}. 
  Assumptions in Lemma~\ref{lemma:conv} are satisfied. Indeed, at each node, the
  number of cost-based cuts, $h_c = \{c^\top \varz \geq \sigma^{[i]}(t) \}$,
  that can be generated along the algorithm evolution is finite, because
  $\sigma^{[i]}(t) = \lceil c^\top \varz^{[i]}(t) \rceil$ and $c^\top
  \varz^{[i]}(0) \leq c^\top
  \varz^{[i]}(t) \leq J^\star$ for all $t$, with $\costval^\star$ being the
  (integer-valued) optimal cost  to~\eqref{eq:MILP}. 
  Therefore, by Lemma~\ref{lemma:conv}, the
  sequences $\{\costval^{[i]}(t)\}_{t \geq 0}$ and
  $\{\varz^{[i]}(t)\}_{t \geq 0}$ at each node $i$ converge in a finite number
  of communication rounds. 
  Now, since the communication graph is jointly strongly connected, by
  Lemma~\ref{lemma:agreement} the sequences converge in a finite number of
  communication rounds to a common limit cost value $\bar{\costval}$ and a
  common limit point $\bar{\varz}$ respectively.

  To prove that $\bar{\costval} = \costval^\star$, we start by highlighting two
  facts.  First, $\costval^{[i]}(t) \leq \costval^\star$, $\forall t \geq 0$,
  and for all $i \in \{ 1, \dots, \Nag\}$, because every agent minimizes the
  cost function over $\Htmp(t)$, which by construction (it is the intersection of
  the collected bases and the local generated cuts) always contains
  $\convS$. Thus, $\bar{\costval} \leq \costval^\star$.
  Second, there exists $\bar{t}$ such that $\varz^{[i]}(t) = \bar{\varz}$ for
  all $t \geq \bar{t}$ and $\forall i \in \{ 1, \dots, \Nag\}$. This means that
  $\bar{\varz}$ satisfies the local constraints $h^{[i]}$ for all
  $i\in\{1, \ldots, \Nag\}$ and, therefore,
  $\bar{\varz} \in \left( \bigcap\limits_{i=1}^{\Nag} h^{[i]}\right) = P$.
  Moreover, since $\bar{z} = (\bar{x}, \bar{y})\in P$ is such that
  $\bar{x} \in \integer^{d_Z}$, it holds $\bar{z} \in \convS$ by
  construction. Thus,
  $\costval(\bar{\varz}) \geq \min_{\varz \in \convS} \costval(\varz) =
  \costval^\star$, giving the cost optimality.
  Finally, being $\bar{\varz}$ feasible and cost optimal, it is the unique lex-optimal
  solution (i.e., $\bar{\varz} = \varz^\star$).%

\subsection{Proof of Theorem~\ref{thm:eps_conv}}
In this subsection, we prove Theorem~\ref{thm:eps_conv}, i.e., $\epsilon$-\DiMILP/ computes an $\epsilon$-suboptimal solution of \ac{MILP}~\eqref{eq:MILP}. 
Before proving the convergence result, we need the following lemma.

\begin{lemma}[Approximate solution] \label{lemma:subepsilon} Suppose
  MILP~\eqref{eq:MILP} satisfies Assumption~\ref{ass:BoundFeas}. Then for any
  given $\epsilon>0$ an optimal solution $(\rho_I^\epsilon, z^\epsilon)$ to
  MILP~\eqref{eq:MILP_epigraphrhoI} exists and $\varz^\epsilon$ is an
  $\epsilon$-suboptimal (feasible) solution to MILP~\eqref{eq:MILP}.
\end{lemma}
\begin{IEEEproof} 
  By Assumption~\ref{ass:BoundFeas}, there exists $z^\star \in \convS$ such that
  $c^\top z^\star \leq c^\top z, \forall z \in \convS$.  If
  $c^\top z^\star/\epsilon \in \integer$, then
  $(\rho_I^\epsilon, z^\epsilon) = (c^\top z^\star/\epsilon, z^\star)$ is a
  feasible, optimal solution of
  \eqref{eq:MILP_epigraphrhoI}. 
  Any optimal solution $(\rho_I^\epsilon, z^\epsilon)$
  of~\eqref{eq:MILP_epigraphrhoI} has cost $c^\top z^\star / \epsilon$ and, by
  construction, $z^\epsilon \in \convS$. Thus,
  $c^\top z^\epsilon / \epsilon \leq c^\top z^\star / \epsilon =
  \rho_I^\epsilon$
  and, being $z^\epsilon \in \convS$, $c^\top z^\epsilon \geq c^\top z^\star$,
  so that $z^\epsilon$ is also optimal for \eqref{eq:MILP}.
  Now suppose $c^\top z^\star/\epsilon \not\in \integer$.  We first show that
  the pair $\rho_I^\epsilon = \lceil c^\top \varz^\star / \epsilon \rceil$,
  $\varz^\epsilon = \varz^\star$ is an optimal solution of
  \eqref{eq:MILP_epigraphrhoI}. Indeed, $\varz^\star$ is feasible for
  \eqref{eq:MILP_epigraphrhoI}, i.e., $\varz^\star \in \convS$, and
  $c^\top \varz^\star \leq \epsilon \lceil c^\top \varz^\star / \epsilon \rceil
  = \epsilon \rho_I^\epsilon$
  (simply because
  $c^\top \varz^\star / \epsilon \leq \lceil c^\top \varz^\star / \epsilon
  \rceil$). 
  Moreover, there does not exist any $\varz \in \convS$ such that
  $c^\top \varz = \epsilon \lfloor c^\top \varz^\star / \epsilon \rfloor$. In
  fact, by the principle of optimality, we have
  $c^\top \varz \geq c^\top \varz^\star $, $\forall z \in \convS$, and
  $c^\top \varz^\star > \epsilon \lfloor c^\top \varz^\star / \epsilon \rfloor$
  (we have already handled the case of integer-valued optimal cost). 
  This proves the existence of an optimal solution of
  \ac{MILP}~\eqref{eq:MILP_epigraphrhoI}.
  Now we
  simply notice that any $(\rho_I^\epsilon, \varz^\epsilon)$ optimal solution of
  \eqref{eq:MILP_epigraphrhoI} is such that $\varz^\epsilon \in \convS$ and
  $c^\top \varz^\epsilon \leq \epsilon \rho_I^\epsilon = \epsilon \lceil c^\top
  \varz^\star / \epsilon \rceil$.
  Moreover, $c^\top \varz^\epsilon \geq c^\top \varz^\star$ because
  $\varz^\epsilon \in \convS$, so that
  $c^\top \varz^\star \leq c^\top \varz^\epsilon \leq \epsilon \lceil c^\top
  \varz^\star / \epsilon \rceil$.  By subtracting $c^\top \varz^\star$, we have
 \[
 0 \leq c^\top \varz^\epsilon - c^\top \varz^\star \leq \epsilon \lceil c^\top
 \varz^\star / \epsilon \rceil - c^\top \varz^\star <\epsilon\,,
 \]
thus concluding the proof.
\end{IEEEproof}

We are now ready to prove that $\epsilon$-\DiMILP/ converges in a finite number
of communication rounds to an $\epsilon$-suboptimal solution of
\ac{MILP}~\eqref{eq:MILP}.

\begin{IEEEproof}[Proof of Theorem~\ref{thm:eps_conv}]
  First, note that we can prove the theorem by replacing the feasible set of
  \ac{MILP}~\eqref{eq:MILP_epigraphrhoI} with
  $P^\epsilon = \{ (\rho_I, \varz) \in \real^{d+1} : a_i^\top \varz \leq b_i, i
  = 1,\ldots,n, \, c^\top \varz \leq \epsilon \rho_I, \rho_I \leq M_\rho\}$,
  where $M_\rho \geq \lceil \max_{\varz\in \PP} c^\top\varz / \epsilon \rceil$.
  Here we are slightly abusing notation, since we denote by $P^\epsilon$ the
  bounded version of the polyhedron. Using $P^\epsilon$ it can be shown that the
  assumptions of Theorem~\ref{thm:int_conv} hold for
  \ac{MILP}~\eqref{eq:MILP_epigraphrhoI}.  Moreover, the optimal cost of
  \ac{MILP}~\eqref{eq:MILP_epigraphrhoI} is integer (indeed, the cost depends
  only on the integer variable $\rho_I$) and the communication graph is jointly
  strongly connected.
  Thus, the sequences 
  $\{ \costval^{[i]}(t)\}_{t \geq 0}$ and
  $\{( \rho_I^{[i]}(t), \varz^{[i]}(t))\}_{t \geq 0}\, i \in \{1,
  \ldots, \Nag \}$, converge, respectively, to the optimal cost, $\costval^\epsilon$, and the 
  lex-optimal solution, $(\rho_I^\epsilon, \varz^\epsilon)$, of MILP~\eqref{eq:MILP_epigraphrhoI}
  in a finite number of communication rounds. 
  Finally, by Lemma~\ref{lemma:subepsilon}, the solution
  ${\varz}^\epsilon$ is an $\epsilon$-suboptimal solution of MILP~\eqref{eq:MILP}.
\end{IEEEproof}

\section{Numerical Computations} %
\label{sec:NumComp}
To corroborate the theoretical analysis and better highlight some appealing
features of our algorithms, we provide numerical computations. 
We concentrate on $\epsilon$-\DiMILP/ since it does not require integer-valued
optimal cost.
First, we consider randomly generated MILP instances and we perform a numerical
Monte Carlo analysis of the algorithm convergence while varying network size
and tolerance $\epsilon$.
Second, we consider a multi-agent multi-task assignment problem that is relevant
for cooperative robotics. For this scenario we test our algorithm in an
unreliable network due to packet loss.

\subsection{Scalability and Convergence Rate for Randomly Generated MILPs}
\begin{figure}[t!]
	\centering
	\subfloat[]{\resizebox{0.48\columnwidth}{!}{%
\includegraphics{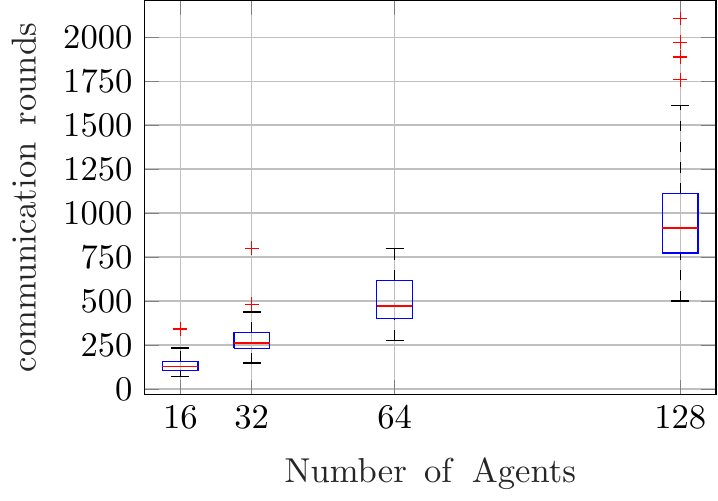}
}\label{fig:boxplot_conv1}}
	\subfloat[]{\resizebox{0.48\columnwidth}{!}{%
\includegraphics{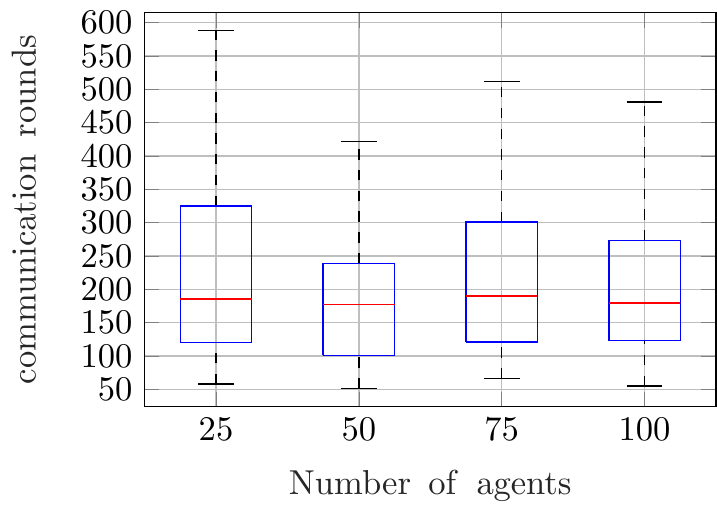}
}\label{fig:boxplot_conv2}}
	\caption{Communication rounds for convergence while varying the
            number of agents for (a) a cyclic digraph with $\diam=\Nag-1$ and
            (b) an Erd\H{o}s-R\'{e}nyi random graph with diameter $\diam=8$.
            Each box plot shows the minimum and maximum communication rounds
            (whiskers), 25\% and 75\% percentiles (lower and upper
            box edges), median (red line), and outliers (red crosses).}
	\label{fig:boxplot_conv}
\end{figure}

We first perform a numerical Monte Carlo analysis to study the number of
communication rounds for convergence while varying the network size.
For each network size, we randomly generate $50$ MILP instances with fixed
dimension, $d=10$ and $d_Z = 3$, and run $\epsilon$-\DiMILP/ with
$\epsilon = 0.1$ and a bounding box $H_M$ with $M = 100$.
The inequality constraints of the MILP instances are randomly generated as follows. 
The vectors $a_i$ are drawn from the standard Gaussian distribution. 
The term $b_i$ is uniformly randomly generated in $[0, 50]$. 
The cost vector $c$ is obtained as $c = a_i^\top \hat{c}$ where $\hat{c}$ is uniformly randomly generated in $[0, 1]^n$. 
The problems generated according to this model are feasible.
The LP relaxation of such a MILP is one of the models proposed
in~\cite{todd1991probabilistic}. 
We first choose a cyclic digraph for which the diameter, $\diam$, is
proportional to the number of agents, specifically $\diam=\Nag-1$, and
consider the following cases: number of agents equal to $16$, $32$, $64$ and
$128$.
The results are shown in Figure~\ref{fig:boxplot_conv1}. 
The red central line of each box shows the median value of the communication
rounds for the $50$ random MILPs. The lower and upper edges of the blue box
represent the $25$-th and $75$-th percentiles.
Some outliers\footnote{Outliers are evaluated as follows. Let $Q_1$ and $Q_3$ be the $25$-th and $75$-th percentile of the samples, respectively. 
A sample is an outlier if it is greater than $Q_3+1.5(Q_3-Q_1)$ or less than $Q_1-1.5(Q_3-Q_1)$.} (red crosses) can be observed. This is not surprising considered
that MILPs are $\mathcal{NP}$-hard.
We point out that the number of communication rounds needed for convergence
grows linearly with the network size (i.e., number of agents and diameter).
Then, we consider Erd\H{o}s-R\'{e}nyi static digraphs with diameter
fixed $\diam=8$, and a growing number of nodes (and so agents), namely $25$, $50$, 
$75$ and $100$.
The results are shown in Figure~\ref{fig:boxplot_conv2}. 
Here we highlight that the median value of communication rounds needed for 
convergence is between $177$ and $190$, that is the completion time of the
algorithm is roughly constant with respect to the number of agents.

\begin{figure}
\centering
\includegraphics{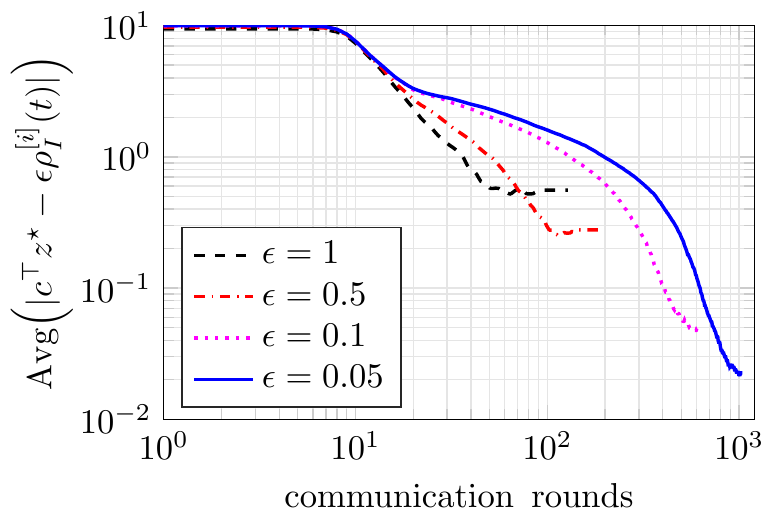}
	\caption{
	Convergence behavior of $\epsilon$-\DiMILP/, 
	with different $\epsilon$ values.
	For each communication round $t$, the figure shows the
	value of $|c^\top\varz^\star-\epsilon\rho^{[i]}_I(t)|$ 
	averaged over all the agents and all the random instances.
	}
	\label{fig:epsilon_conv_rate}
\end{figure}

Second, we perform a numerical Monte Carlo analysis to study the behavior of
  $\epsilon$-\DiMILP/ while varying $\epsilon$. We randomly generate $50$
  \ac{MILP} instances with $d=10$, $d_Z = 3$ and $256$ constraints, as
  before. The constraints are distributed among $\Nag =64$ agents, so that each
  one knows $4$ constraints. Agents communicate according to an
  Erd\H{o}s-R\'{e}nyi random graph with $\diam=7$.  Each instance is solved with
  $\epsilon=0.05,0.1,0.5,1$, and the results are shown in
  Figure~\ref{fig:epsilon_conv_rate}.
For each communication round $t$, we plot the mean value (over instances and
  agents), $\mathrm{Avg}(|c^\top\varz^\star-\epsilon\rho^{[i]}_I(t)|)$, 
  of the distance
  $|c^\top\varz^\star-\epsilon\rho^{[i]}_I(t)|$ of each instance from the optimal cost.
  Different curves are associated to different
  values of $\epsilon$.
  As expected, by decreasing $\epsilon$, the cost at convergence is closer to
  the optimal one, but a larger number of communication rounds is needed for
  convergence.
To speed-up the convergence, a ``multi-resolution'' strategy may be implemented,
in which $\epsilon$ is reduced at each run and the previous computed solution is
used to initialize the new run.  %
Notice that, if agents can apply one of the distributed stopping criteria
discussed in the previous section, e.g., if the communication graph is static,
this scheme can be implemented in a purely distributed way.
Finally, we point out that the algorithm exhibits a convergence rate behavior
similar to centralized cutting-plane algorithms. That is, in the first
iterations the cost soars to the optimal cost value. As the algorithm evolves,
one can observe the typical tailing-off effect characterizing cutting planes
algorithms, which slows down the convergence rate.
How to combine cutting planes with other schemes, to avoid this tailing-off, is
definitely an interesting research direction.

\subsection{Distributed Multi-Agent Multi-Task Assignment under Unreliable Network }
In this section we apply $\epsilon$-\DiMILP/ to a Multi-Agent Multi-Task
  Assignment problem for robotic networks.  Here we use the variant of
  $\epsilon$-\DiMILP/ described in Remark~\ref{rmk:multi_cut} in which multiple
  cuts are generated.  

We consider a group of heterogeneous mobile robots that need to accomplish a
set of tasks.
Each task must be performed at a given target location and each vehicle has the
capability to execute only some of the tasks.
Robots move in a constrained space including obstacles, which (for
  simplicity) are modeled as rectangles.
Given the capability constraints and forbidden areas, the goal is to assign the tasks among the
vehicles in order to minimize the mission completion time which is defined as the
time needed for the last vehicle to finish its mission, see,
e.g.,~\cite{bellingham2003multi}. 
We assume that valid paths, from vehicles to target locations, are available
  together with the corresponding finishing times.  This could model, e.g., a
  structured environment, as in factories or warehouses, in which vehicles
  follow pre-defined paths on the floor of the operating area.  Thus, we setup
  the following multi-agent multi-task assignment problem.

Given $N_\theta$ target locations,
$\{T_1, \ldots, T_{N_\theta} \}$, $N_\nu$ vehicles,
$\{V_1, \ldots, V_{N_\nu} \}$, and $\Npaths$ paths, we formulate the
multi-agent multi-task assignment problem as the following MILP:
\begin{align}\label{eq:TaskAssignment}
\min_{\varint, \vareal} &\; y &&  \nonumber\\
\subj &\; \vecthetaT \, \varint \leq - w_\theta \,, && \forall \theta \in \{ 1,\ldots, N_\theta \}  \nonumber \\
&\;  \vecnuT \varint \leq 1 \,,\, && \forall \nu \in \{ 1,\ldots, N_\nu  \} \\
&\; \vecgammaT \varint \leq \vareal \,, && \forall \gamma \in \{ 1,\ldots, \Npaths \} \nonumber 
\end{align}%
where $\varint \in \integer^{\Npaths}$ is a vector of binary decision variables
associated to the paths ($\varint_\indexa$ is equal to one if the $\indexa$-th path is
selected, and $0$ otherwise), while $\vareal \in \real$ is a continuous decision
variable associated to the mission completion time.
The vector $\vectheta \in \real^{\Npaths}$ is a vector whose $\indexa$-th entry
is $1$ if target location $T_\theta$ is visited by path $\indexa$ and $0$
otherwise.
The vector $\vecnu \in \real^{\Npaths}$ is a vector whose $\indexa$-th entry is
$1$ if path $\indexa$ is assigned to vehicle $V_\nu$ and $0$ otherwise.
Each vector $\vecgamma \in \real^{\Npaths}$ has the form
$\vecgammaT = [0, \ldots, 0, \tau_\indexa, 0, \ldots, 0]$, where $\tau_\indexa$
is the completion time for the $\indexa$-th path.
The first set of constraints enforces that each target location is visited at
least $w_\theta$ times, the second one prevents more than one path being
assigned to each vehicle, and the third one forces $y$ to be the overall
finishing time, i.e., $\vareal = \max_{\indexa} \tau_\indexa$. %
It is worth noting that an optimal solution of~\eqref{eq:TaskAssignment} cannot
be obtained as the solution of its LP relaxation, i.e., $\varint \in \real^{\Npaths}$
and $0 \leq \varint_\indexa \leq 1$, $\forall \indexa \in \{1, \ldots, \Npaths \}$,
due to the presence of the capability constraints.

We consider $\Nag$ agents, $\{A_1, \ldots, A_{\Nag} \}$, randomly located in
  the operating region and communicating according to a jointly strongly
  connected communication graph.
The graph is built by considering a proximity criterion, i.e., two agents
  are connected if their distance is less than a threshold.

\begin{figure}[ttpb]
  \begin{center}
	
	\subfloat[]{%
	\includegraphics[scale=.75]{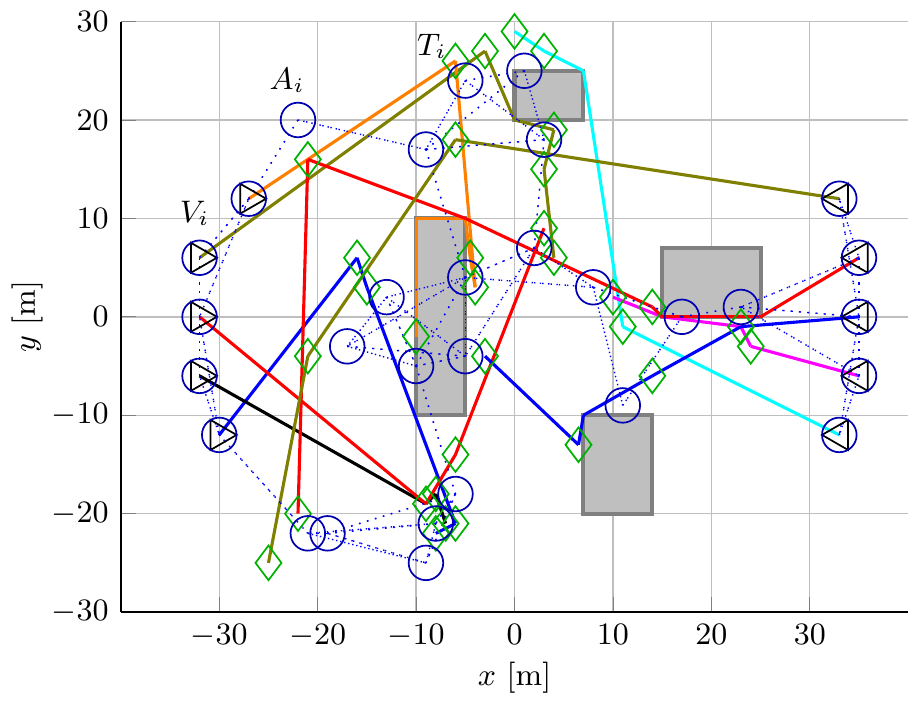}%
	\label{fig:TA1}} 
	
	\subfloat[]{\includegraphics[scale=.4]%
	{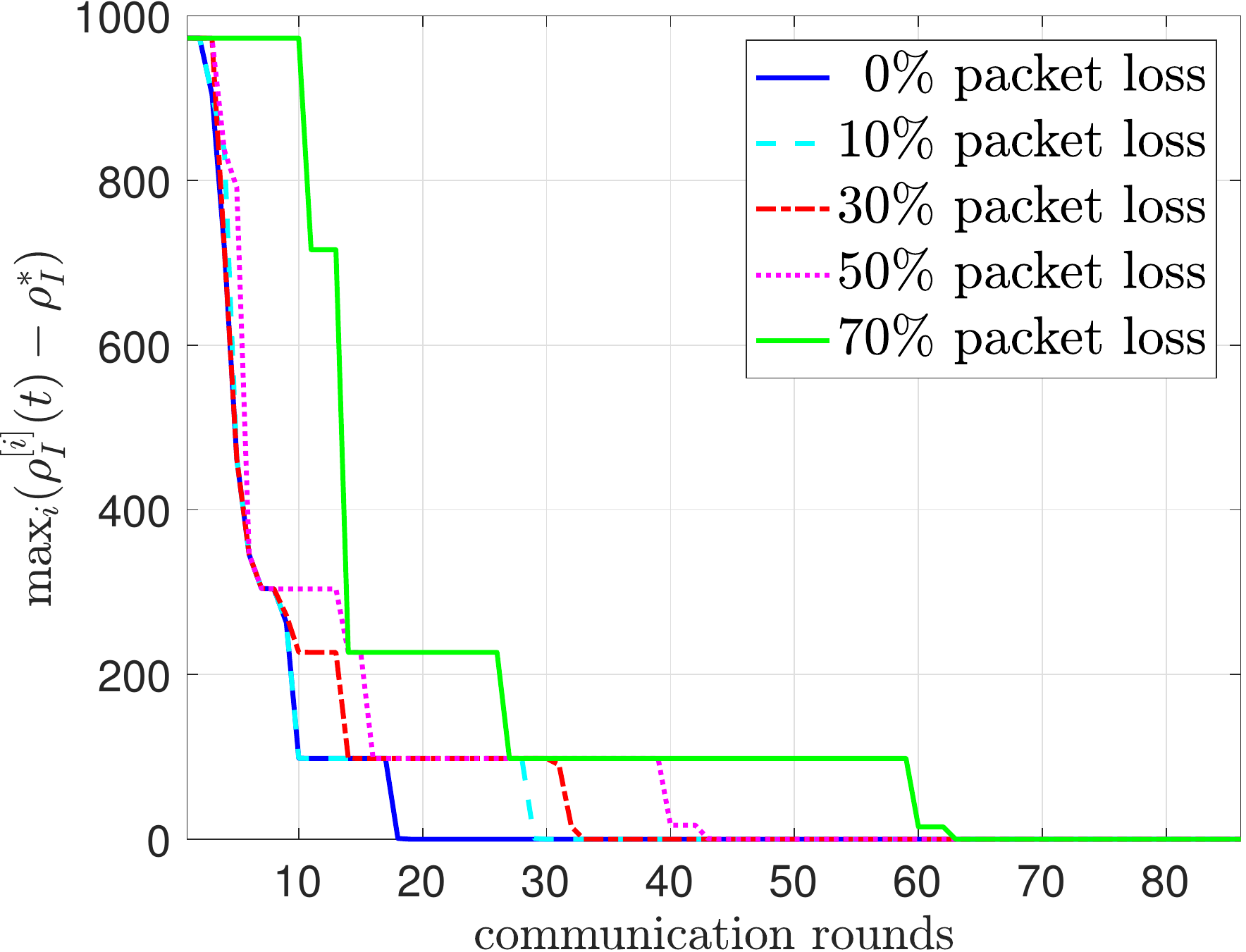}\label{fig:TA2}} 	
		\caption{In (a) Multi-Agent Multi-Task assignment solution: $10$ vehicles
          (black triangles), $32$ tasks (green diamonds), $10$ optimal paths (colored lines)
          and  $30$ agents (blue circles) connected by a proximity graph (blue dotted lines).
         In (b)
         $\max_{i \in \{1,\ldots, \Nag\}} (\rho^{[i]}_I(t)-\rho^\star_I$) for
         different percentages of packet loss. 
	}
	\label{fig:TA}
  \end{center}
\end{figure}

Each agent only knows the constraints in problem~\eqref{eq:TaskAssignment}
associated to paths traversing a neighboring area. 
For example, in Figure~\ref{fig:TA1}, agent $A_i$ (blue circle) only knows the
paths traversing a circle area of radius $14$m centered at its
position, and communicates with other agents in this circle.
In this framework, a robotic vehicle, e.g., $V_i$ in Figure~\ref{fig:TA1},
  is also an agent participating to the computation.

In order to show the robustness of $\epsilon$-\DiMILP/ against failures in the
communication network, we introduce a probability of packet loss for each edge
in the fixed communication graph.
Specifically, at each communication round, an agent $i$ does not receive a
message from a neighboring agent $j$ if a random variable in $[0,1]$ is smaller
than a fixed threshold. 
The threshold defines the percentage of packet loss for
edge $(j,i)$.
We consider $0\%$ (ideal case, no packet loss), $10\%$, $30\%$, $50\%$ and
$70\%$ of packet loss. 
For each packet loss percentage, we run $\epsilon$-\DiMILP/ on an instance of
problem~\eqref{eq:TaskAssignment} with $\epsilon = 0.1$, $N_\theta = 32$,
  $N_\nu = 10$, $\Npaths = 71$, and $\Nag = 30$. The solution is depicted in
  Figure~\ref{fig:TA1} (colored lines).
  The convergence results are shown in Figure~\ref{fig:TA2} in which we plot the
  difference between $\max_{i \in \{1,\ldots, \Nag\}} \rho^{[i]}_I(t)$ and the
  solution $\rho^\star_I$ of MILP~\eqref{eq:MILP_epigraphrhoI} (computed by
  using a solver in YALMIP~\cite{lofberg2004yalmip}).
  For the ideal case (no packet loss), agents reach consensus on the optimal
  cost after $19$ %
  communication
  rounds. %
When increasing the percentage of packet loss, as expected, the number of
communication rounds increases, see Figure~\ref{fig:TA2}. %
Notably, as from the theory, the algorithm converges also for very high
  percentages of packet loss as, e.g., $70\%$.

\section{Conclusion}
\label{sec:concl}
In this paper, we proposed an algorithmic framework, with two specific
algorithms, to solve Mixed-Integer Linear Programs over networks.
In the proposed distributed setup, the constraints of the \ac{MILP} are assigned to a
network of agents. The agents have a limited amount of memory and computation
capabilities and are able to communicate with neighboring agents.
Following the idea of centralized cutting-plane methods for \ac{MILPs}, each
agent solves local LP relaxations of the global problem, generates cutting
planes, and exchanges active constraints (a candidate basis) with neighbors. 
We proved that agents running the first algorithm, \texttt{INT}-\DiMILP/, reach
consensus on the lex-optimal solution of the \ac{MILP}, under the assumption of
integer-valued optimal cost, in a finite number of communication rounds.
Removing the assumption of integer-valued optimal cost, we proved that 
agents running the algorithm $\epsilon$-\DiMILP/ reach consensus on 
a suboptimal solution (up to a given tolerance $\epsilon$) of the \ac{MILP}, 
in a finite number of communication rounds. 
Both algorithms involve low-cost local computations and work under asynchronous,
unreliable, and directed networks.
Finally, we performed a set of numerical computations suggesting that the
completion time of $\epsilon$-\DiMILP/ scales nicely with the
diameter of the communication graph. We also tested the algorithm on a
multi-agent multi-task assignment setup in unreliable networks.
  Future investigations may include the solution of \ac{MILPs} with a large
  number of decision variables and the combination of pure cutting plane methods
  with other tools, e.g. branch and bound, as in centralized schemes.

\section*{Acknowledgment}
The authors would like to thank Andrea Camisa for the deep discussions and useful suggestions.

\begin{small}
% Generated by IEEEtran.bst, version: 1.14 (2015/08/26)

\end{small}
\begin{IEEEbiography}
  [{\includegraphics[width=2.5cm]{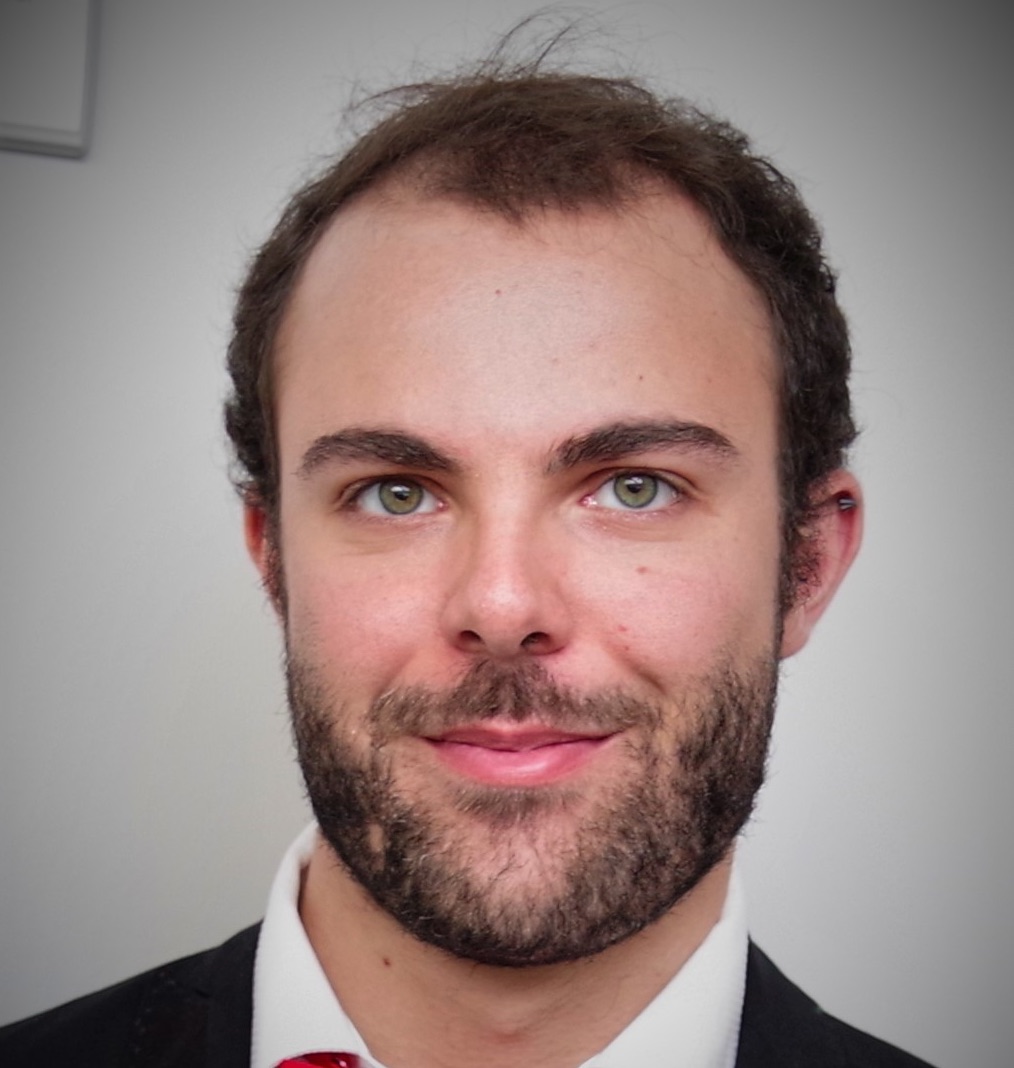}}]{Andrea Testa}
received the Laurea degree ``summa cum laude'' in Computer Engineering from the Universit\`a del Salento, Lecce, Italy in 2016.

He is a Ph.D. student in Engineering of Complex Systems at the Universit\`a del Salento, Lecce, Italy since November 2016. He was a visiting scholar at LAAS-CNRS, Toulouse, (July to September 2015 and February 2016). 

His research interests include control of UAVs and distributed optimization.
\end{IEEEbiography}
\vspace{-3.5cm}
\begin{IEEEbiography}
  [{\includegraphics[width=2.5cm]{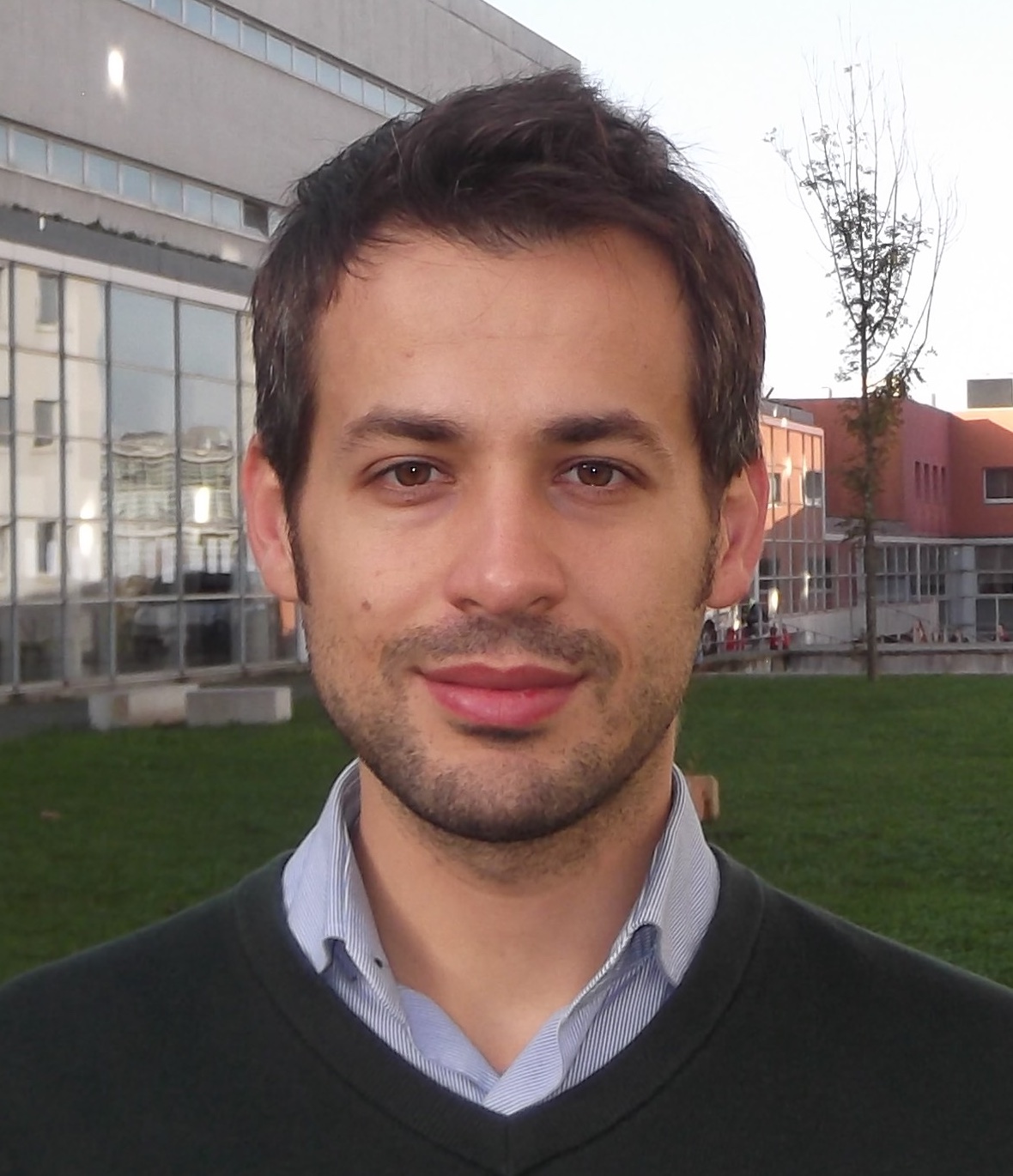}}]{Alessandro Rucco}
received the Master's degree in automation engineering and the Ph.D. degree in information engineering from 
the Universit\`a del Salento, Lecce, Italy, in 2009 and 2014, respectively. 

He was Visiting Student at the 
University of Paris Sud Sup\'elec (March-April 2010) and at the University of Colorado Boulder 
(January-September 2012). He was Team Leader of the VI-RTUS team winning the International Student
Competition Virtual Formula 2012. 
From August 2014 to September 2016, he was a Post-Doctoral Researcher at the Department of Electrical 
and Computer Engineering, University of Porto, Portugal. From October 2016 to July 2018, he was a Post-Doctoral Researcher at the 
Department of Engineering, Universit\`a del Salento. 

His research interests include distributed optimization,  
nonlinear optimal control, modeling, trajectory optimization, and high-performance maneuvering for car and 
aerial vehicles. 
\end{IEEEbiography}
\vspace{-3.5cm}
\begin{IEEEbiography}
  [{\includegraphics[width=2.5cm]{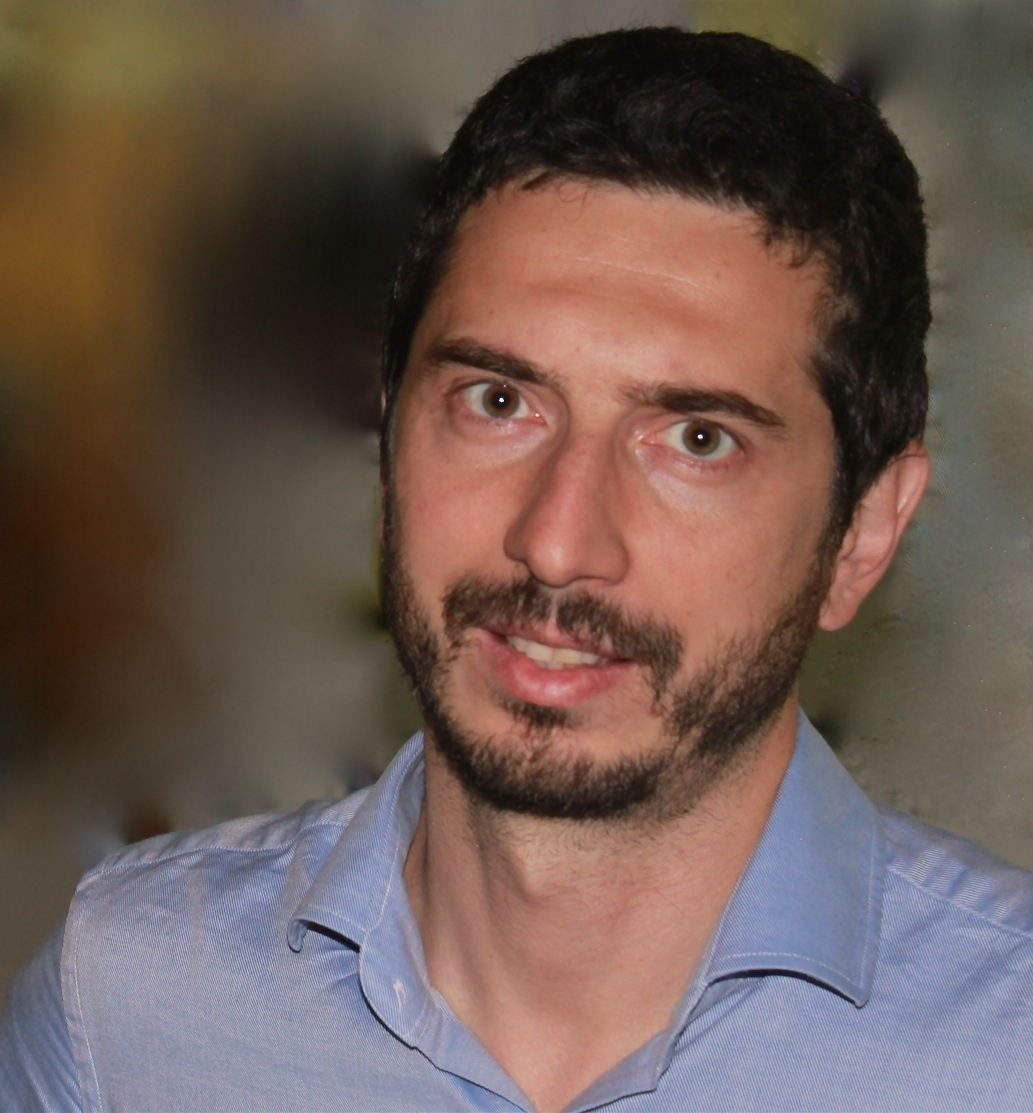}}] {Giuseppe Notarstefano} 
received the Laurea degree summa cum laude in electronics engineering from the Universit\`a di Pisa, Pisa, Italy, in 2003 and the Ph.D. degree in automation and operation research from the Universit\`a di Padova, Padua, Italy, in 2007.

He is a Professor with the Department of Electrical, Electronic, and Information Engineering G. Marconi, Alma Mater Studiorum Universit\`a di Bologna, Bologna, Italy. He was Associate Professor (from June 2016 to June 2018) and
previously Assistant Professor, Ricercatore (from February 2007), with the Universit\`a del Salento, Lecce, Italy. He has been Visiting Scholar at the University of Stuttgart, University of California Santa Barbara, Santa Barbara, CA, USA and University of Colorado Boulder, Boulder, CO, USA. His research interests include distributed optimization, cooperative control in complex networks, applied nonlinear optimal control, and trajectory optimization and maneuvering of aerial and car vehicles.

Dr. Notarstefano serves as an Associate Editor for \emph{IEEE Transactions on Automatic Control}, \emph{IEEE Transactions on Control Systems Technology}, and \emph{IEEE Control Systems Letters}. He is also part of the Conference Editorial Board of IEEE Control Systems Society and EUCA. He is recipient of an ERC Starting Grant 2014.
\end{IEEEbiography}

\end{document}